\documentstyle[amssymb,10pt]{amsart}

\newtheorem{prop}{Proposition}[section]
\newtheorem{prop:def}{Proposition-Definition}[section]
\newtheorem{defin}{Definition}[section]
\newtheorem{lemma}{Lemma}[section]
\newtheorem{thm}{Theorem}[section]

\theoremstyle{remark}
\newtheorem{remark}{Remark}

\newtheorem{example}{Example}

\begin{document}

\newcommand{\nc}{\newcommand} \nc{\on}{\operatorname}

\nc{\¦}{{|}}

\nc{\pa}{\partial} \nc{\cA}{{\cal A}} \nc{\cB}{{\cal B}} \nc{\cC}{{\cal C}}
\nc{\cE}{{\cal E}} \nc{\cF}{{\cal F}} \nc{\cG}{{\cal G}} \nc{\cH}{{\cal H}} 
\nc{\cK}{{\cal K}} \nc{\cL}{{\cal L}} \nc{\cM}{{\cal M}} \nc{\cN}{{\cal N}}
\nc{\cO}{{\cal O}} \nc{\cP}{{\cal P}} \nc{\cR}{{\cal R}} \nc{\cU}{{\cal U}} 
\nc{\cV}{{\cal V}} \nc{\cW}{{\cal W}} \nc{\cX}{{\cal X}} \nc{\cD}{{\cal D}}

\nc{\sh}{\on{sh}} \nc{\Id}{\on{Id}} \nc{\Diff}{\on{Diff}}
\nc{\id}{\on{id}}

\nc{\ul}{\underline}\nc{\ol}{\overline}

\nc{\ad}{\on{ad}}
\nc{\Der}{\on{Der}}
\nc{\End}{\on{End}}
\nc{\Lie}{\on{Lie}}
\nc{\res}{\on{res}}
\nc{\ddiv}{\on{div}}
\nc{\FS}{\on{FS}}\nc{\Hitch}{{\on{Hitch}}}
\nc{\GNR}{{\on{GNR}}}\nc{\AJ}{{\on{AJ}}}\nc{\BNR}{{\on{BNR}}}

\nc{\card}{\on{card}} \nc{\dimm}{\on{dim}}
 \nc{\tr}{\on{tr}}
 \nc{\gr}{\on{gr}}
\nc{\Jac}{\on{Jac}} \nc{\Ker}{\on{Ker}} \nc{\Den}{\on{Den}} \nc{\genus}{\on{genus}}
 \nc{\Imm}{\on{Im}}
\nc{\limm}{\on{lim}} \nc{\Ad}{\on{Ad}}

\nc{\ev}{\on{ev}}
\nc{\Hol}{\on{Hol}}
\nc{\Det}{\on{Det}}

\nc{\Cone}{\on{Cone}}
\nc{\pseudo}{{\on{pseudo}}}

\nc{\class}{{\on{class}}}
\nc{\rat}{{\on{rat}}}
\nc{\local}{{\on{local}}}
\nc{\an}{{\on{an}}}

\nc{\Lift}{{\on{Lift}}}
\nc{\Mer}{{\on{Mer}}}
\nc{\mer}{{\on{mer}}}

\nc{\lift}{{\on{lift}}}
\nc{\diff}{{\on{diff}}}
\nc{\Aut}{{\on{Aut}}}

\nc{\DO}{{\on{DO}}}
\nc{\Frac}{{\on{Frac}}}
\nc{\cl}{{\on{class}}}

\nc{\Fil}{{\on{Fil}}}

\nc{\Bun}{\on{Bun}}
\nc{\diag}{\on{diag}}
\nc{\KZ}{{\on{KZ}}}

\nc{\CB}{{\on{CB}}}
\nc{\out}{{\on{out}}}
\nc{\Hom}{{\on{Hom}}}

\nc{\FO}{{\on{FO}}}

\nc{\al}{\alpha}
\nc{\de}{\delta}
\nc{\si}{\sigma}
\nc{\ve}{\varepsilon}
\nc{\z}{\zeta}

\nc{\vp}{\varphi}
\nc{\la}{{\lambda}}
\nc{\g}{\gamma}
\nc{\eps}{\epsilon}

\nc{\PsiDO}{\Psi\on{DO}}\nc{\BM}{{\on{BM}}}

\nc{\AAA}{{\mathbb A}}\nc{\CC}{{\mathbb C}}\nc{\KK}{{\mathbb K}}
\nc{\NN}{{\mathbb N}}\nc{\PP}{{\mathbb P}}\nc{\RR}{{\mathbb R}}
\nc{\VV}{{\mathbb V}}\nc{\ZZ}{{\mathbb Z}}\nc{\EE}{{\mathbb E}}
\nc{\TT}{{\mathbb T}}

\nc{\bla}{{\mathbf \lambda}}
\nc{\bv}{{\mathbf v}}

\nc{\bz}{{\mathbf z}}
\nc{\bt}{{\mathbf t}}

\nc{\bP}{{\mathbf P}}
\nc{\kk}{{\mathbf k}}
\nc{\A}{{\mathfrak a}}
\nc{\B}{{\mathfrak b}}
\nc{\G}{{\mathfrak g}}
\nc{\HH}{{\mathfrak h}}
\nc{\mm}{{\mathfrak m}}
\nc{\N}{{\mathfrak n}}

\nc{\SG}{{\mathfrak S}}
\nc{\La}{\Lambda}

\nc{\wt}{\widetilde}

\nc{\wh}{\widehat}
\nc{\bn}{\begin{equation}}
\nc{\en}{\end{equation}}

\nc{\SL}{{\mathfrak{sl}}}
\nc{\ttt}{{\mathfrak{t}}}
\nc{\GL}{{\mathfrak{gl}}}

\title[Quantizations of Hitchin and Beauville-Mukai integrable systems]
{Quantizations of the Hitchin and Beauville-Mukai integrable systems}

\begin{abstract}
Spectral transformation is known to set up a birational morphism between 
the Hitchin and Beauville-Mukai integrable systems. The corresponding phase spaces 
are: (a) the cotangent bundle of the moduli space of 
bundles over a curve $C$, and (b) a symmetric power of the cotangent surface $T^*(C)$. 
We conjecture that this morphism can be quantized, and we check 
this conjecture in the case where $C$ is a rational curve with marked points
and rank $2$ bundles. We discuss the relation of the resulting isomorphism of
quantized algebras with Sklyanin's separation of variables. 
\end{abstract}

\author{B. Enriquez and V. Rubtsov}

\address{B.E.: IRMA (CNRS), Universit\'e  Louis Pasteur, 7, rue Ren\'e Descartes,
67084 Strasbourg, France}

\address{V.R.: D\'epartement de Math\'ematiques, Universit\'e d'Angers,
49045 Angers, France}

\maketitle

\subsection*{Introduction} 

In \cite{Hitchin}, Hitchin introduced an algebraically completely integrable 
system (ACIS). This system is attached to a complex curve $C$ and a reductive 
complex Lie group $G$. Its phase space is the cotangent bundle to the 
moduli space of stable $G$-bundles over $C$. This system can 
be generalized to the case of marked points and parabolic bundles.
It is important because its quantization (in a $\cD$-module sense)
plays a central role in the Beilinson-Drinfeld study of the 
geometric Langlands conjectures. The Hitchin system is also a common
framework for several many-body integrable systems. 

The fact that the Hitchin system is algebraically completely integrable
means that its fibers are open subsets of suitable Jacobians. This fact
was proved by Hitchin in \cite{Hitchin} and relies on a variant of the 
spectral construction. This construction attaches to each point of 
the moduli space, a covering $C_h$ of $C$ (contained in the surface $T^*C$) 
and a line bundle $\cL$ over $C_h$. 

In \cite{GNR}, Gorsky, Nekrasov and the second-named author studied
the case $G = \on{GL}_n(\CC)$. They constructed a bijection between 
open subsets of (a) the restriction of 
the Hitchin system to the connected component $\cU_{n,d}$ of bundles of degree 
$d = n^2(g-1)+1$, and (b) a Beauville-Mukai system. The latter system 
is an integrable system associated to a Poisson surface; its phase space 
is a symmetric power of the surface (see \cite{Cortona}, 
\cite{Bottacin}, \cite{Mukai}). 

The isomorphism of \cite{GNR} is constructed as follows: 

(1) Hitchin's spectral transformation is a map from $\TT^*\cU_{n,d}$ 
to the moduli space of pairs $(C_h,\cL)$, 

(2) Hurtubise constructs a birational morphism from this moduli 
space to the $\wt g = d$-th symmetric power of $T^*(C)$ (\cite{Hurt}). 

In this paper, we first come back to the construction of \cite{GNR}. 
We prove that the Hitchin system for bundles of degree $\wt g$ and the 
Beauville-Mukai system are birationally equivalent (Section \ref{section:1}).
We then discuss quantizations of both systems (in the "classical" sense, 
meaning a commutative subalgebra of a noncommutative algebra). A quantization 
of the Beauville-Mukai system was obtained in \cite{EO,ER} (see also \cite{BT} for
another proof of the first results of \cite{ER}). On the other hand, 
the Beilinson-Drinfeld construction gives a quantization of the 
Hitchin system. We conjecture that the isomorphism between both 
classical systems can also be quantized. 
 
In the next sections, we check this conjecture explicitly in the 
particular case where $C$ is $\PP^1$ with marked points and $n = 2$. 
We recall the explicit form of the Hitchin and Beauville-Mukai
systems, as well as of the isomorphism $\GNR$ in Section \ref{sect:1}. 
$\GNR$ is then an isomorphism of Poisson fields with $\CC^\times$-actions.
These fields are the fraction fields of graded Poisson algebras, 
which correspond to "ringed Lie algebras", which we introduce in Section 
\ref{sect:ringed}. We show that $\GNR$ induces an isomorphism of 
ringed Lie algebras (Section \ref{sect:iso:ringed}) and of Poisson 
fields (Section \ref{sect:5}).   
In Section \ref{sect:quant}, we study the quantization of these systems. 
The quantization of the Hitchin system is provided by the Gaudin 
hamiltonians (see \cite{FFR}), and we quantize the Beauville-Mukai 
system according to \cite{ER}. We then obtain the quantization of 
$\GNR$ in Theorem \ref{thm:quant:GNR}. 

Our computations are closely related with Sklyanin's 
separation of variables. We discuss the relations with the works of 
Sklyanin and E.\ Frenkel in Section \ref{sect:skl}. 

Sklyanin's separation of variables was generalized in \cite{Gekht} 
to the case of the group $\on{SL}_n(\CC)$ (at the classical level). 
This work might help to generalize the "explicit" part of our paper
to $n>2$.


\section{A birational equivalence of ACIS} \label{section:1}

In this section, we recall the construction of \cite{GNR} of an isomorphism 
between the Hitchin and Beauville-Mukai systems in the case of a general 
curve $C$. We then 
propose a conjecture on the quantization of this isomorphism. 

\subsection{An isomorphism of ACIS and a quantization conjecture}

Let $C$ be a smooth, connected, complete algebraic curve over $\CC$. 
We denote by $K$ its canonical bundle. When $E$ is a bundle over 
a variety $V$, we denote by $\EE$ the total space of $E$. 

\medskip \noindent
{\it (a) The Hitchin system.} Let $n,d$ be integers such that $n\geq 1$. 
Let $\cU_{n,d}$ be the moduli space of stable bundles over $C$, 
with rank $n$ and degree $d$. Let us set 
$$
\cH = \bigoplus_{i=1}^n H^0(C,K^{\otimes i}). 
$$
The Hitchin fibration is the morphism 
$$
\Hitch : \TT^* \cU_{n,d} \to \cH, 
$$
defined at the level of points by $\Hitch(E,\phi) 
= (\tr(\wedge^i(\phi)))_{i = 1,\ldots,n}$. Here $E$ is a stable 
bundle over $C$ and $\phi\in H^0(C,\End(E)\otimes K)$. $\Hitch$
is a $\CC^\times$-equivariant, Lagrangian fibration (\cite{Hitchin}). 

\medskip \noindent
{\it (b) The Beauville-Mukai system.} Let $S$ be a Poisson surface, let 
$\wt g$ be an integer $\geq 1$, and let $\cL$ be a line bundle over 
$S$, such that $h^0(\cL) = \wt g +1$. Then there is a morphism 
$$
\BM : S^{[\wt g]} \supset U \to \PP(H^0(S,\cL)), 
$$
taking the class of a generic collection $(s_1,\ldots,s_{\wt g})$
of points of $S$ to the line of all sections of $\cL$ vanishing at these 
points. When $S^{[\wt g]}$ is equipped with the Mukai-Bottacin Poisson 
structure (\cite{Mukai,Bottacin}), $\BM$ is a Lagrangian fibration. 

\medskip 

Let us construct $S$ as follows. Let $\overline\KK$ be the compactification 
$\PP(\cO_C \oplus K)$ of $\KK$, where the curve at infinity is blown 
down to a point: so $\overline\KK$ is a one-point compactification of 
$\KK$. Let $C_\infty$ the the zero-section of $\KK \subset \overline\KK$. 
We take $S = \overline\KK$, $\cL = \cO_{\ol\KK}(n C_\infty)$, $\wt g = n^2(g-1)+1$. 
Then 
$$
H^0(\overline\KK,\cO_{\ol\KK}(nC_\infty)) = \bigoplus_{i=0}^n H^0(\overline\KK,
\cO_{\ol\KK}(nC_\infty))^i = \bigoplus_{i=0}^n H^0(C,K^{\otimes i}), 
$$
where the exponent $i$ denotes the subspace of regular functions, 
homogeneous of degree $i$ in each fiber. So we get a dominant morphism 
$$
\BM : \ol\KK^{[\wt g]} \supset U \to \PP(\CC \oplus \cH). 
$$ 
$\cH$ is an open subset of $\PP(\CC \oplus \cH) = \cH \cup \PP(\cH)$, 
so $\BM$ restricts to a morphism 
$$
\BM : \ol\KK^{[\wt g]} \supset U' \to \cH,  
$$ 
where $U'$ is a Zariski-open subset of $\ol\KK^{[\wt g]}$. 
$\ol\KK - C_\infty$ is a cone over $C$ and $\cH$ is graded, so 
both sides of this map are equipped with an action of $\CC^\times$; then 
$\BM$ is $\CC^\times$-equivariant.

\medskip \noindent
{\it (c) The birational morphism $\GNR$.}

\begin{thm} (see \cite{GNR}) \label{GNR}
In $n,d$ are related by $d = n^2(g-1) +1$ and if $\wt g = d$, then there is a 
birational Poisson morphism 
$$
\GNR : \ol\KK^{[\wt g]} \supset U'' \to \TT^* \cU_{n,d}, 
$$
such that the diagram 
$$
\begin{array}{ccccccc}
\ol\KK^{[\wt g]} & \supset & U'' &  & \stackrel{\GNR}{\to} &   & \TT^*\cU_{n,d} \\
   &  &  & \scriptstyle{\BM}\searrow  & & \scriptstyle{\Hitch}\swarrow  & \\
                 &          &    &          &  \cH    &          &                 
\end{array}
$$
commutes. 
\end{thm}

In other words, $\GNR$ sets up a birational morphism of ACIS between the 
Hitchin and Beauville-Mukai systems. At the level of sets of points, 
$\GNR$ and its inverse can be described as follows. 

\medskip 

Let us first describe $\GNR$. Let $(s_1,\ldots,s_{\wt g})$ be a $\wt g$-uple
of generic points of $\KK$. Let $h$ be its image under $\BM: \ol\KK^{[\wt g]} 
\supset U' \to \cH$. Let $C_h$ be the curve of $\¦nC_\infty\¦$ corresponding to $h$; 
then $C_h$ contains the points $(s_1,\ldots,s_{\wt g})$ (and is uniquely
determined by this requirement). When $(s_1,\ldots,s_{\wt g})$ are generic, 
$C_h$ is smooth. Let $p_h : C_h \to C$ be the canonical projection. Its degree
is $n$. Let $E = (p_h)_*(\cO_{C_h}(s_1 + \cdots + s_{\wt g}))$, then $E$ is a rank 
$n$ bundle over $C$. Moreover, we have a natural element $\phi \in 
H^0(C,\End(E) \otimes K)$, provided by the embedding $C_h \subset \KK$. 

\medskip Let us now describe the map $\GNR^{-1}$. Let $(E,\phi)\in 
\TT^*\cU_{n,d}$; let $h$ be the image of $(E,\phi)$ be the Hitchin 
map $\TT^*\cU_{n,d}\to \cH$. When $(E,\phi)$ are generic, the curve 
$C_h$ is smooth. The spectral line bundle $\cL$ over $C_h$ is defined 
as the kernel of $p_h^*(\phi) - \la : p_h^*(E) \to p_h^*(E\otimes K)$, 
where $p_h : C_h \to C$ is the canonical projection, and $\la$ is the 
tautological section of the bundle $p^*(K)$ ($p : \KK \to C$ 
is the canonical projection). 
When $(E,\phi)$ is generic, $\cL$ has degree $\wt g$, and $h^0(C_h,\cL) = 1$. 
The divisor of a nonzero section $\sigma$ of $\cL$ is then an element of 
$C_h^{(\wt g)}$. We define an element of $\KK^{[\wt g]}$ as the ideal of 
length $\wt g$, preimage of the ideal sheaf of $\sigma$ by $\cO_\KK 
\twoheadrightarrow \cO_{C_h}$; this is $\GNR^{-1}(E,\phi)$.

\medskip \noindent
{\it (d) The quantization conjecture.} The algebraic translation of the 
statements of Section (a) is that the rational function fields 
$R(\ol\KK^{[\wt g]})$ and $R(\TT^*\cU_{n,d})$ are Poisson, and we have linear maps 
$\cH^* \to R(\ol\KK^{[\wt g]})$ and $\cH^* \to R(\TT^*\cU_{n,d})$
with Poisson-commutative images. The translation of Theorem \ref{GNR}
is that there is an isomorphism of Poisson fields $R(\TT^*\cU_{n,d})
\to R(\ol\KK^{[\wt g]})$, such that the diagram 
\begin{equation} \label{comm:diag:fields}
\begin{array}{ccccc}
R(\TT^*\cU_{n,d}) & \;         & \to & \;        & R(\ol\KK^{[\wt g]}) \\
\;                & \nwarrow   & \;  & \nearrow  & \;  \\
\;                & \;         & \cH & \;        & \;      
\end{array}
\end{equation}
commutes. In Section \ref{quant:conj}, we will present a conjecture on quantization 
of this diagram. We then check this conjecture in the case when $(n,d) = (2,-1)$
and $C = \PP^1$ with marked points. 

\subsection{Proof of Theorem \ref{GNR}: construction of the birational 
morphism $\GNR$}

\subsubsection{Moduli spaces} We recall some facts about moduli spaces of bundles
(\cite{Sesh,Nitsure}). Recall first the notion of a coarse moduli space. 
Let $F$ be a functor from the category of $\CC$-schemes to that of  sets, a 
{\it coarse moduli space} for $F$ is the data of (1) a $\CC$-scheme $\cM_F$, 
(2) a natural transformation $F \to \on{Mor}_{\on{Schemes}}(-, \cM_F)$, 
inducing a bijection $F(\on{Spec}(\CC)) = \cM_F(\CC)$, and such that for any scheme
$S$, any natural transformation $F \to \on{Mor}_{\on{Schemes}}(-,S)$ gives 
rise to as schemes morphism $\cM_F\to S$, such that the resulting diagram 
$$
\begin{array}{ccc}
F  & \to        & \on{Mor}_{\on{Schemes}}(-, \cM_F) \\
\; & \searrow   & \downarrow                           \\
\; & \;         & \on{Mor}_{\on{Schemes}}(-, S) 
\end{array}
$$ 
commutes. 

Let $n,d$ be integers, with $n\geq 1$. Let $F_{n,d}$ be the functor, defined 
at the level of objects by 
$F_{n,d}(X) = \{$bundles $\cE$ over $X\times C$, such that for any point $x\in X(\CC)$, 
$\cE_{\¦ \{x\}\times C}$ is stable, of rank $n$ and degree $d\} / \sim$, where 
the equivalence relation $\sim$ is defined by $\cE \sim \cE'$ iff there exists a 
line bundle $\cL$ over $X$, such that $\cE$ and $\cE' \otimes p_1^*(\cL)$ are isomorphic
($p_1,p_2$ are the natural projections $X\times C\to X$, $X\times C \to C$). 

\begin{thm}(\cite{Sesh}) There exists a coarse moduli space $\cU_{n,d}$
for $F_{n,d}$. It is a smooth, connected quasiprojective variety
of dimension $n^2(g-1)+1$. 
\end{thm}

Let $G_{n,d}$ the functor defined at the level of objects by 
$G_{n,d}(X) = \{$stable pairs $(\cE,\phi)$ of (a) a bundle $\cE$ over $X\times C$, 
such that for any $x\in X(\CC)$, $\cE_{\{x\}\times C}$ has rank $n$ and degree $d$, 
(b) a bundle morphism $\phi : \cE \to \cE \otimes p_2^*(K)\} / \sim'$, where 

$\bullet$ "stable pair" means that for any $\phi$-invariant proper subbundle $\cF$ 
of $\cE$, we have $\mu(\cF) < \mu(\cE)$ (where $\mu =$ degree/rank). So if $\cE$ 
is stable, any pair $(\cE,\phi)$ is stable; 

$\bullet$ the equivalence relation $\sim'$ is defined 
by $(\cE,\phi) \sim' (\cF,\psi)$ iff there exists a line bundle 
$\cL$ over $X$, and an isomorphism of $(\cF,\psi)$ with $(\cE\otimes 
p_2^*(\cL),\phi \otimes \on{id}_{p_2^*(\cL)})$. 

\begin{thm} \label{thm:openness} (\cite{Nitsure})  
There exists a moduli space $\cM_{n,d}$
for $G_{n,d}$. It is a smooth, quasiprojective variety, which contains 
$\TT^*\cU_{n,d}$ as an open subvariety. 
\end{thm}

One checks that $\TT^*\cU_{n,d}(\CC)$ identifies with the preimage 
$\pi^{-1}(\cU_{n,d}(\CC))$ by the (set-theoretic) projection 
$\cM_{n,d}(\CC)\to \{$classes of bundles over $C\}$, taking $(E,\phi)$
to $E$. 

Recall now some openness results: 

\begin{thm} (see \cite{Sesh,Nitsure}) Let $X$ be a complex variety. 

1) Let $\cE$ be a bundle over $X\times C$, then the subset of $X(\CC)$
of all points $x$ such that $\cE_{\¦\{x\}\times C}$ is stable, is Zariski-open. 

2) Let $(\cE,\phi)$ be a pair of a bundle over $X\times C$ and a bundle morphism 
$\phi : \cE \to \cE \otimes p_2^*(K)$, then the subset of $X(\CC)$ of all 
points $x$, such that $(\cE,\phi)_{\¦\{x\}\times C}$ is stable, is 
Zariski-open. 
\end{thm}

\subsubsection{Construction of a bundle $\cE$ over $U''' \subset \ol\KK^{[\wt g]}$}

The Beauville-Mukai map is a dominant morphism 
$$
\BM : \ol\KK^{[\wt g]} \supset U' \to \cH, 
$$
where $U'\subset\ol\KK^{[\wt g]}$ is Zariski-open, and $\cH = \oplus_{i=1}^n 
H^0(C,K^{\otimes i})$. 

To each $h\in \cH$, one attaches the spectral curve 
$C_h$. Recall this construction. Let  $U_C\subset C$ be an 
open subset, over which $K$ is trivial, and let $\omega : 
\cO_{C\¦U_C} \to K_{\¦U_C}$ be a trivialization. Then $\omega$
induces an identification $p^{-1}(U_C) \to U_C \times \AAA^1$. We define 
$\varphi_{U_C} \in \cO_{\cH} \otimes \cO_{p^{-1}(U_C)}$ as the image of 
the function $\wt\varphi \in \cO_{\cH} \otimes \cO_{C\times \AAA^1}$, 
which takes $((h_1,\ldots,h_n),x,y)$ to 
$$
\wt\varphi = y^n + (h_1/\omega)(x) y^{n-1} + \cdots + (h_n/\omega^n)(x). 
$$
For $h$ fixed in $\cH$, the ideals $(\varphi_{U_C}(h,-))$ glue to an ideal sheaf 
defining $C_h$. When $C_h$ is smooth, it has genus $\wt g$. 
When $P_1,\ldots,P_{\wt g}\in \KK$ are distinct and $(P_1,\ldots,P_{\wt g})\in U'$, 
then $C_h$ contains the points $P_1,\ldots,P_{\wt g}$. So 
$$
m_{i,\wt g +1} \big( \varphi_{U_C} \circ (\BM,\on{id}) \big) = 0,  
$$
where $m_{i,\wt g +1} : \cO_{\KK^{\wt g +1}} \to \cO_{\KK^{\wt g}}$ 
is the multiplication of $i$th and $(\wt g+1)$th factors. 

There exists a Zariski-open subset $\cH_{\on{smooth}} \subset \cH$, such that 
for any $h\in \cH_{\on{smooth}}$, the curve $C_h$ is smooth. 
According to \cite{BNR}, remark 3.5, $\cH_{\on{smooth}}$ is nonempty
if $K^n$ admits a section without multiple zeroes; this is always the case
according to \cite{Hitchin}, 5.1 (based on the Bertini theorem). Then 
$U_{\on{smooth}} = \BM^{-1}(\cH_{\on{smooth}})$ is a Zariski-open subset of $U'$. 

Let $\pi : \KK^{[\wt g]} \to C^{(\wt g)}$ be the canonical projection, 
let $\Delta \subset C^{(\wt g)}$ be the image of the union of all diagonals
under $C^{\wt g} \to C^{\wt g} / \SG_{\wt g} = C^{(\wt g)}$, and let $\cV = 
\pi^{-1}(C^{(\wt g)} - \Delta)$. Set $U''' = \cV \cap U_{\on{smooth}}$. 
Then $U'''$ is an open subset of $\KK^{[\wt g]}$. 

We construct a sheaf $\cE$ of $\cO_{U'''\times C}$-modules over $U''' \times C$
as follows. Let $U_1,\ldots,U_{\wt g}$ be open subsets of $\KK$, and let $U_C$ 
be an open subset of $C$ over which $K$ trivializes. Then 
$$
O(U_1,\ldots,U_{\wt g}) = \pi'\Big( \big(\cup_{\sigma\in\SG_{\wt g}}
U_{\sigma(1)} \times \cdots \times U_{\sigma(\wt g)}\big) 
\cap (\pi')^{-1}(U''')\Big) 
$$
is an open subset of $U'''$. 
Here $\pi'$ is the natural projection $\KK^{\wt g} - \on{diagonals}
\to \KK^{[\wt g]}$. 
We set 
\begin{align*}
& \Gamma(O(U_1,\ldots,U_{\wt g}) \times U_C,\cE) = 
\{\on{regular\ functions\ }f\on{\ on\ } 
\\& \Big( \big(\cup_{\sigma\in\SG_{\wt g}}
U_{\sigma(1)} \times \cdots \times U_{\sigma(\wt g)}
 \big) \cap (\pi')^{-1}(U''')\Big) 
\times p ^{-1}(U_C) \subset \KK^{\wt g +1}, 
\\ & 
\on{symmetric\ in\ their\ }\wt g\on{\ first\ arguments,\ such\ that\ for\ any}
\\ & i = 1,\ldots,\wt g,\ m_{i,\wt g +1}(f) = 0 \} / \on{ideal\ generated\ by\ }
\varphi_{U_C}\circ (\BM,\on{id}). 
\end{align*}

\begin{prop} \label{prop:local:freeness}
$\cE$ is a locally free sheaf on $U''' \times C$ of rank $n$. 
If $(P_1,\ldots,P_{\wt g})\in U'''$, then its restriction to 
$\{(P_1,\ldots,P_{\wt g})\}\times C$ is isomorphic to 
$(p_h)_*(\cO(P_1 + \cdots + P_{\wt g}))$, where 
$h = \BM(P_1,\ldots,P_{\wt g})$.  
\end{prop}

{\em Proof.} The stalk of $\cE$ at $((P_1,\ldots,P_{\wt g}),P)\in U''' \times C$ is
$$
\cE_{((P_1,\ldots,P_{\wt g}),P)} = \limm_{\rightarrow} \Gamma(V,\cE), 
$$
where $V$ runs over all open subsets containing $((P_1,\ldots,P_{\wt g}),P)$. 
Set 
$$
\wh\cE_{((P_1,\ldots,P_{\wt g}),P)}  = 
\cE_{((P_1,\ldots,P_{\wt g}),P)}  \otimes_{\cO_{((P_1,\ldots,P_{\wt g}),P)}}
\wh\cO_{((P_1,\ldots,P_{\wt g}),P)} , 
$$
where $\cO_{((P_1,\ldots,P_{\wt g}),P)}$ is the local ring at 
$((P_1,\ldots,P_{\wt g}),P)$, and $\wh \cO_{((P_1,\ldots,P_{\wt g}),P)}$ is 
its completion. 
According to \cite{Hartshorne}, ex.\ 5.7, $\cE$ is locally free iff 
each $\cE_{((P_1,\ldots,P_{\wt g}),P)}$ is free, and since 
$\cE_{((P_1,\ldots,P_{\wt g}),P)}$ is finitely generated, 
according to \cite{Hartshorne}, Lemma 8.9, it is free of rank $n$ iff 
$\wh\cE_{((P_1,\ldots,P_{\wt g}),P)}$ is. 
So we should prove that 
$\wh\cE_{((P_1,\ldots,P_{\wt g}),P)}$ is a free 
$\wh\cO_{((P_1,\ldots,P_{\wt g}),P)}$-module of rank $n$. 

Let us compute  $\wh\cE_{((P_1,\ldots,P_{\wt g}),P)}$. 
We identify locally $\KK$ with a product $C\times \AAA^1$, 
and derive from there local coordinates $(z_1,\la_1),\ldots,(z_{\wt g},\la_{\wt g})$
of $\KK$ at $P_1,\ldots,P_{\wt g}$. Let $z$ be a formal coordinate of $C$ at $P$. 

If $U_P$ is a formal neighborhood of $P$ at $C$, then $p^{-1}(U_P)$ identifies
with a trivial bundle $U_P \times \AAA^1$ with coordinates $(z,\la)$. Then the function 
$\wt \varphi$ takes the form 
$$
\wt\varphi(\la,z\¦ z_1,\la_1,\ldots,z_{\wt g},\la_{\wt g}) = 
\la^n + \sum_{i=1}^n \la^{n-i} \al_i(z\¦ z_1,\la_1,\ldots,z_{\wt g},\la_{\wt g}),  
$$
where each $\al_i$ belongs to $\CC[[z,z_1,\ldots,\la_{\wt g}]]$. 

The projections $P'_1,\ldots,P'_{\wt g}$ of $P_1,\ldots,P_{\wt g}$ on $C$
are all distinct. Assume first that none of them is equal to $P$. Then the 
vanishing condition $m_{i,\wt g +1}(f) = 0$ is empty. So 
$\wh\cE_{((P_1,\ldots,P_{\wt g}),P)}$ identifies with the quotient
$$
\CC[\la][[z,z_1,\la_1,\ldots,z_{\wt g},\la_{\wt g}]] / 
(\la^n + \sum_{i=1}^n \la^{n-i} \al_i(z,z_1,\ldots,\la_{\wt g})). 
$$
This is a free $\CC[[z,z_1,\ldots,\la_{\wt g}]]$-module with basis 
$(1,\la,\ldots,\la^{n-1})$. 

Assume now that one projection $P'_i$ equal $P$. Then all $P'_j$, $j\neq i$, 
are different from $P$, and $P'_i = (P,\la_0)$, with $\la_0\in\CC$;  
$z_i$ is a local coordinate on $C$ at $P$ and $\la_i$ is a local coordinate
on $\AAA^1$ at $\la_0$. 
There are unique elements $\beta_i$ of $\CC[[z,z_1,\ldots,\la_{\wt g}]]$, 
such that 
\begin{equation} \label{A}
\wt\varphi(\la,z\¦z_1,\la_1,\ldots,z_{\wt g},\la_{\wt g}) = 
(\la - \la_i)^n + \sum_{\al = 1}^n (\la - \la_i)^{n-\al} \beta_\al
(z,z_1,\ldots,\la_{\wt g}).  
\end{equation}
Moreover, $\wt \varphi(\la_i,z_i \¦ z_1,\la_1,\ldots,z_{\wt g},\la_{\wt g}) = 0$, 
so $\beta_0(z_i,z_1,\ldots,\la_{\wt g}) = 0$. So there exists 
$\gamma\in \CC[[z,z_1,\ldots,\la_{\wt g}]]$, such that 
\begin{equation} \label{B}
\beta_0(z,z_1,\ldots,\la_{\wt g}) = (z-z_i)
\gamma(z,z_1,\ldots,\la_{\wt g}) . 
\end{equation}
Now 
\begin{align} \label{starbis}
& \nonumber 
\{f\in\CC[\la][[z,z_1,\ldots,\la_{\wt g}]] \¦ m_{i,\wt g+1}(f) = 0\}
\\& 
= (\oplus_{\al\geq 1} \CC (\la - \la_i)^\al) [[z,z_1,\ldots,\la_{\wt g}]]
\oplus (z-z_i) \CC[[z,z_1,\ldots,\la_{\wt g}]]. 
\end{align}

\begin{lemma}
The quotient 
$$
\{f\in\CC[\la][[z,z_1,\ldots,\la_{\wt g}]] \¦ m_{i,\wt g+1}(f) = 0\}
/ (\wt \varphi(\la,z\¦ z_1,\ldots,\la_{\wt g}))
$$
is a free module over $\CC[[z,z_1,\ldots,\la_{\wt g}]]$ with basis 
\begin{equation} \label{family}
(z-z_i,\la-\la_i,\ldots,(\la-\la_i)^{n-1}). 
\end{equation}
\end{lemma}

{\em Proof of Lemma.} Let us show that (\ref{family}) is a generating
family. According to (\ref{starbis}), it suffices to show that each 
$(\la-\la_i)^k$, $k\geq 1$ is a combination of the elements of (\ref{family}). 
We show this by induction on $k$: this is obvious when $k\in \{1,\ldots,n-1\}$. 
Let us assume that we have shown that for some $\al_{k,l} \in 
\CC[[z,z_1,\ldots,\la_{\wt g}]]$, we have 
$$
(\la-\la_i)^k = \Big( \sum_{l=1}^{n-1} \al_{k,l}(z,\ldots,\la_{\wt g})
(\la-\la_i)^l \Big) + \al_{k,0}(z,\ldots,\la_{\wt g})(z-z_i). 
$$ 
Then 
\begin{align*}
& (\la-\la_i)^{k+1} = \Big( \sum_{l=2}^{n-1} \al_{k,l-1}(z,\ldots,\la_{\wt g})
(\la-\la_i)^l \Big) + \al_{k,0}(z,\ldots,\la_{\wt g}) (z-z_i)(\la-\la_i)
\\ & 
-\al_{k,n-1}(z,\ldots,\la_{\wt g}) \Big( 
\sum_{\al = 1}^{n-1} (\la-\la_i)^{n-\al} \beta_\al(z,z_1,\ldots,\la_{\wt g})
+(z-z_i) \gamma(z,z_1,\ldots,\la_{\wt g}) \Big) 
\end{align*}
according to (\ref{A}) and (\ref{B}). So (\ref{family}) is generating.  

On the other hand, since the degree in $\la-\la_i$ of a 
combination of the elements of (\ref{family}) is always $<n$, 
it cannot belong to $(\wt\varphi(\la,z\¦z_1,\ldots,\la_{\wt g}))$   
without vanishing. So (\ref{family}) is also a free family. 
\hfill \qed \medskip

This ends the proof of the local freeness of $\cE$. The second statement of
Proposition \ref{prop:local:freeness} is clear. \hfill \qed\medskip 

\subsubsection{A bundle morphism $\phi : \cE \to \cE \otimes p_2^*(K)$}

Recall that $p : \KK \to C$ is the canonical projection. 
We have a section $\sigma$ of $p^*(K)$, defined as follows: 
for $c\in C$, the restriction to $p^{-1}(c)$ of $p^*(K)$ identifies
with the trivial bundle over $T_c(C)^*$ with fiber $T_c(C)^*$, and the 
restriction of $\sigma$ is the identity map $T_c(C)^*\to T_c(C)^*$. 
Then the operation 
$$
f(P_1,\ldots,P_{\wt g},P) \mapsto f(P_1,\ldots,P_{\wt g},P) \sigma(P), 
$$
where $(P_1,\ldots,P_{\wt g},P) \in\KK^{\wt g +1}$, induces a linear map 
$$
\Gamma(O(U_1,\ldots,U_{\wt g})\times U_C,\cE) \to 
\Gamma(O(U_1,\ldots,U_{\wt g})\times U_C,\cE\otimes p_2^*(K)),   
$$
which defines a bundle morphism $\phi : \cE \to \cE \otimes p_2^*(K)$. 

\subsubsection{The map $\GNR : U'' \to \TT^*\cU_{n,d}$}

For $h\in \cH_{\on{smooth}}$, let $(p_h)_* : 
\Jac^{\wt g}(C_h) \to \Bun_n(C)$ be the direct image map. 
Let $J^{\wt g}_h$ be the preimage $(p_h)_*^{-1}(\cU_{n,d})$
of stable bundles. Then the proof of Theorem 1 of \cite{BNR} says that for 
$h$ in a nonempty open of $\cH$, $J^{\wt g}_h \to \cU_{n,d}$ is dominant. Since 
$\dimm \Jac^{\wt g}(C_h) = \dimm \cU_{n,d}$, $J^{\wt g}_h$ contains a 
nonempty open subset of $\Jac^{\wt g}(C_h)$. Since the same is true of the image of 
the Abel-Jacobi map $\AJ_{\wt g} : C_h^{(\wt g)} \to \Jac^{\wt g}(C_h)$, 
$\AJ_{\wt g}(C_h^{(\wt g)})$ and $J_h^{\wt g}$ have a dense intersection. 
So there exists $(P_1,\ldots,P_{\wt g})\in U'''$, such that 
$(p_h)_*(\cO(P_1 + \cdots + P_{\wt g}))$ is stable, where 
$h = \BM(P_1,\ldots,P_{\wt g})$. 

It then follows from Theorem \ref{thm:openness}, 1), that 
$U'' = \{u\in U''' \¦ \cE_{\{u\}\times C}\on{\ is\ stable}\}$
is a nonempty open subset of $U'''$. Now since $\cU_{n,d}$ is a coarse
moduli scheme, the bundle $\cE$ gives rise to a morphism $U'' \to \cU_{n,d}$. 

Moreover, the each $u\in U''$, the pair $(\cE,\phi)_{\¦\{u\}\times C}$ 
is also stable, so according to \cite{Nitsure}, $(\cE,\phi)_{\¦U''\times C}$
gives rise to a morphism $\mu : U'' \to \cM_{n,d}(K)$. According to \cite{Nitsure}, 
the subset of $\cM_{n,d}(K)(\CC)$ of all stable pairs $(\cE,\phi)$ such that 
$\cE$ is stable, can be identified with $\TT^*\cU_{n,d}$, so $\mu$
factors through a morphism $\GNR : U'' \to \TT^*\cU_{n,d}$. 

\subsubsection{$\BM = \Hitch \circ \GNR$}

Let us now show that the diagram 
\begin{equation} \label{GNR:diag}
\begin{array}{ccccc}
U'' &          & \stackrel{\GNR}{\to}& & \TT^*\cU_{n,d} \\
      & \scriptstyle{\BM}{\searrow} &   & {\swarrow}\scriptstyle{\Hitch} & \\
      &          &          \cH       &          & 
\end{array}
\end{equation}
commutes. 
$\BM$ and $\Hitch\circ \GNR$ are two morphisms $U'' \to \cH$
of algebraic varieties. The maps 
$\BM(\CC)$ and $(\Hitch\circ\GNR)(\CC) : U''(\CC) \to \cH(\CC)$
between sets of closed points are the same. Therefore 
$\BM = \Hitch\;\! \circ \;\!\GNR$. 

\subsubsection{$\GNR$ is a birational morphism} \label{sect:Zariski}

In the last sections, we have constructed an open subset $U'' \subset 
\ol\KK^{[\wt g]}$ and a morphism $U'' \to \TT^*\cU_{n,d}$. To prove that 
it induces a birational morphism between $\ol\KK^{[\wt g]}$ and 
$\TT^*\cU_{n,d}$, we use the following lemma (see \cite{Hartshorne}). 

\begin{lemma}
1) Let $f : X \to Y$ be a morphism of connected complex algebraic varieties. 
Assume that there exists a smooth point $x_0\in X(\CC)$, such that $f(x_0)$
is smooth and $df(x_0)$ is a linear isomorphism. 
Then $f$ is dominant and corresponds to a finite extension of 
rational function fields $R(Y) \subset R(X)$. 

2) If in addition, for any $y\in Y(\CC)$, $f^{-1}(y)$ is connected, then $f$
has degree $1$, so $R(Y) = R(X)$, i.e., $f$ is a birational 
morphism. 
\end{lemma}

{\em Proof.} 1) follows from the Zariski main theorem (\cite{Dieud}, 
Theorem 1, p.\ 134). 2) follows from Corollary 1 of {\it loc.\ cit.} 
\hfill \qed \medskip 

Let us apply this lemma to $\GNR : U'' \to \TT^*\cU_{n,d}$. For this, 
we should check the following statements: 

\begin{lemma} \label{lemma:1} There exists $(P_1,\ldots,P_{\wt g})
\in U''(\CC)$, such that $d(\GNR)(P_1,\ldots,P_{\wt g})$ is a linear
isomorphism.  
\end{lemma}

\begin{lemma} \label{lemma:2} For any $(E,\phi) \in \TT^*\cU_{n,d}(\CC)$, 
$\GNR^{-1}(E,\phi)$ is connected. 
\end{lemma}

{\em Proof of Lemma \ref{lemma:1}.} Recall that the diagram 
(\ref{GNR:diag}) commutes. $\BM : \ol\KK^{[\wt g]} \supset U' \to \cH$
is dominant, so the rank of its differential is maximal on some nonempty open 
subset $O \subset U'''$. Let us fix $h\in \BM(O)$. Then the rank of 
$d(\BM)$ is maximal on a nonempty open of $\BM^{-1}(h)\cap U'''$. 
The restriction of $\GNR$ to $\BM^{-1}(h)\cap U'''$ induces a morphism
$$
\GNR' : \BM^{-1}(h) \cap U''' \to \Hitch^{-1}(h). 
$$
Let $C_h\subset \KK$ be the curve of the linear system $\¦nC_\infty\¦$
attached to $h$. Then $\BM^{-1}(h) \cap U'''$ identifies with a 
nonempty open of $C_h^{(\wt g)}$. Moreover, we have a commutative diagram 
\begin{equation} \label{useful:diag}
\begin{array}{ccc}
C_h^{(\wt g)}  &  \stackrel{\AJ_{\wt g}}{\to} & \Jac^{\wt g}(C_h)  \\
 \cup         &                             &  \cup        \\
\BM^{-1}(h) \cap U''' & \to & J^{\wt g}_h \\
 & \scriptstyle{\GNR'} \searrow & \downarrow \scriptstyle{\BNR}\\
& & \Hitch^{-1}(h) 
\end{array}
\end{equation}
where $\AJ_{\wt g}$ is the Abel-Jacobi map, $\Jac^{\wt g}(C_h)$
is the degree $\wt g$ Jacobian of $C_h$, and $\BNR$ is the map defined in 
\cite{BNR}. 

Since $\BNR$ is bijective (see \cite{BNR}, Proposition 3.6), and $\AJ_{\wt g}$
is a birational morphism, $\GNR'$ is a birational morphism. So we can 
find $u\in \BM^{-1}(h)\cap U'''$, such that 
$$
d(\GNR')(u) : T_u(\BM^{-1}(h') \cap O) \to T_{\GNR'(u)}(\Hitch^{-1}(h)) 
$$ 
is a linear isomorphism. We have a commutative diagram 
$$
\begin{array}{ccccccccc}
0 \to & T_u(\BM^{-1}(h)\cap O) & \to & T_u(O) & \stackrel{d(\BM)(u)}{\to} 
& T_h(\cH) & \to & 0 \\
 & \downarrow &  &  \scriptstyle{d(\GNR)(u)}\downarrow &  &  
\scriptstyle{\id}\downarrow &  &  \\
0 \to & T_{\GNR(u)}(\Hitch^{-1}(h)) & \to & T_{\GNR(u)}(\TT^*\cU_{n,d}) 
& \to & T_h(\cH) & \to & 0 
\end{array}
$$
where the rows are exact. Since $d(\BM)(u)$ is surjective, $d(\GNR)(u)$ is a
linear isomorphism. \hfill \qed \medskip 

{\em Proof of Lemma \ref{lemma:2}.} Let us recall the structure of the 
Abel-Jacobi map in degree $\wt g$ (see \cite{Fay}). The morphism 
$$
\AJ_{\wt g} : C_h^{(\wt g)} \to \Jac^{\wt g}(C_h)
$$
is surjective; it is 1-1 on all of $\Jac^{\wt g}(C_h)$, except on the 
codimension $2$ divisor equal to $\on{div}(K) - \AJ_{\wt g - 2}(C_h^{(\wt g - 2)})$, 
where $\on{div}(K)$ is the canonical divisor class of $C_h$ and $\AJ_{\wt g-2}$
is the degree $\wt g-2$ Abel-Jacobi map. The fiber of $\AJ_{\wt g}$ at a point  
$\on{div}(K) - \AJ_{\wt g -2}(x_1,\ldots,x_{\wt g -2})$ of the codimension $2$
divisor is connected. Indeed, it is the image of the map 
$$
\PP(V_{x_1,\ldots,x_{\wt g -2}}) \to C_h^{(\wt g)}, 
$$
where $V_{x_1,\ldots,x_{\wt g-2}} = \{$regular $1$-differentials on 
$C_h$, vanishing at $x_1,\ldots,x_{\wt g-2}\}$, taking $\omega$ to 
$\on{div}(\omega) - (x_1 + \cdots + x_{\wt g-2})$ (here $\on{div}(\omega)$
is the divisor class of $\omega$). 

Let us now fix $(E,\phi) \in \TT^*\cU_{n,d}$ and let us study 
$\GNR^{-1}(E,\phi)$. If this set is nonempty, then $\Hitch(E,\phi)\in
\cH_{\on{smooth}}$. Then $\GNR^{-1}(E,\phi)$ identifies with the fiber
$(\GNR')^{-1}(E,\phi)$ under
$$
\GNR' : \BM^{-1}(h)\cap U'' \to \Hitch^{-1}(h). 
$$
It follows from the diagram (\ref{useful:diag}), from the fact that
$\BNR$ is bijective, from the structure of $\AJ_{\wt g}$ recalled
above, and from the fact that a Zariski open subvariety of a connected
variety is connected, that $(\GNR')^{-1}(E,\phi)$ is connected. 

This proves Lemma \ref{lemma:2}. \hfill \qed \medskip 

\subsection{Quantization}

In this section, we construct quantizations of the rational function fields
$R(\TT^* \cU_{n,d})$ and $R(\ol\KK^{[\wt g]})$. For this, we first construct
quantizations of the underlying graded rings, and then quantize their
fraction fields using ``algebraic microlocalization''\! 
. We then propose a conjecture on quantization of the isomorphism of Theorem 
\ref{GNR}. 

\subsubsection{Fraction fields of graded commutative rings}  

Let $A = \oplus_{i\in \ZZ} A_i$ be a graded commutative algebra, 
equipped with a Poisson structure of degree $-1$. We assume that $A$ is
integral. Let $Q(A)$ be the fraction field of $A$. It contains the 
graded subalgebra $Q(A)^{\on{graded}}$, which can either be viewed as 
the localization of $A$ with respect to the multiplicative part
$(\cup_{i\in\ZZ} A_i) - \{0\}$, or as the direct sum $\oplus_{i\in\ZZ}
Q(A)_i$, where $Q(A)_i$ is the subspace of $Q(A)$ of all 
homogeneous of degree $i$ elements. 

Let us denote by $\wh{Q(A)}$ the completion of $Q(A)^{\on{graded}}$
with respect to the family $(\oplus_{j\leq i} Q(A)_j)_{i\in\ZZ}$. 
So an element of $\wh{Q(A)}$ is a sequence $(q_i)_{i\in \ZZ}$, 
such that $q_i\in Q(A)_i$ and $q_i = 0$ for $i$ large enough. 
We have a double inclusion 
$$
Q(A)^{\on{graded}} \subset Q(A) \subset \wh{Q(A)}; 
$$
the two first algebras are dense in the last one; all three
algebras are Poisson, and the Poisson structure of $Q(A)^{\on{graded}}$
has degree $-1$; only the two last algebras are fields. 

\begin{example} \label{ex:1}
Set $A = \{$rational functions on $\TT^*\cU_{n,d}$, polynomial in the fibers$\}$, 
then we have $A = \oplus_{i\geq 0}A_i$, where $A_i = \{$elements of $A$, 
homogeneous of degree $i$ in the fibers$\}$. Here $A_i = 0$ for $i<0$. 
Then for $i\in\ZZ$, $Q(A)_i = \{$rational functions on $\TT^*\cU_{n,d}$, 
rational and homogeneous of degree $i$ in the fibers$\}$.   
\end{example}

\begin{example} \label{ex:2}
Set $A' = \oplus_{i\geq 0} H^0(C,K^{\otimes i})$, $A'_{-i} = 
H^0(C,K^{\otimes i})$ for $i\geq 0$ and $A'_i = 0$ for $i>0$. 
Then $A' = \oplus_{i\in\ZZ} A'_i$ is a graded algebra, with Poisson structure of
degree $-1$ (see \cite{EO}). Let us set $A'' = S^{\wt g}(A') = 
((A')^{\otimes \wt g})^{\SG_{\wt g}}$, then $A''$ is also graded, 
and its Poisson structure also has degree $-1$. The elements of $A'$
are the regular functions on $\ol\KK - C_\infty$, so they are polynomial 
in the fibers of the projection $\KK - C_\infty \to C$. 
The elements of $Q(A'')_i$ are the rational 
functions on $\KK^{[\wt g]}$, rational and homogeneous of degree $i$
in the fibers of $(\KK - C_\infty)^{[\wt g]} \to C^{(\wt g)}$.  
\end{example}

\subsubsection{Quantization of fraction fields}

By a quantization of the graded Poisson algebra $A = \oplus_{i\in\ZZ} A_i$, 
we understand a filtered algebra $B$, with filtration 
$$
\cdots \subset B_i \subset B_{i+1} \subset \cdots, 
$$ 
with $\cap_{i\in\ZZ} B_i = \{0\}$, $\cup_{i\in\ZZ} B_i = B$, 
$B_i B_j \subset B_{i+j}$, $B = \limm_{\leftarrow} (B/B_i)$, 
i.e., $B$ is complete for the topology defined by the $(B_i)_{i\in\ZZ}$, 
such that the associated graded $\gr(B)$ of $B$ is commutative, 
and we have an isomorphism of graded Poisson algebras $\gr(B)
\stackrel{\sim}{\to} A$.  

The theory of ``algebraic microlocalization'' provides 
a quantization of $Q(A)^{\on{graded}}$ starting from a quantization of $A$
(see \cite{Spr}).

Elements of $Q(B)$ are infinite series $\sum_{i\geq 0} P_i Q_i^{-1}$, 
where $\omega(P_i) - \omega(Q_i) \to - \infty$. The product of two such series
is computed using the formula
$$
Q^{-1} R = R Q^{-1} + [R,Q] Q^{-2} + \cdots + \ad(-Q)^{n-1}(R)Q^{-n} + \cdots
$$
Such a series is zero iff its degree is $\leq d$ for any $d$. 
To prove that a series $\sum_{i\geq 0} P_i Q_i^{-1}$ 
is zero, one should iterate the following 
operation: (1) show that the image of the series in $Q(A)_\omega$
is zero, where $\omega = \on{max}_i(\omega(P_i) - \omega(Q_i))$, 
and (2) using this vanishing, rewrite it as a series $\sum_{i\geq 0} P'_i (Q'_i)^{-1}$, 
where  $\on{max}_i(\omega(P'_i) - \omega(Q'_i)) < \on{max}_i(\omega(P_i) - \omega(Q_i))$.

\subsubsection{Quantization of the Hitchin and Beauville-Mukai phase
spaces}

A quantization of the algebra $A$ of Example \ref{ex:1} is the algebra 
$\Diff(\cU_{n,d})$ of rational differential operators on $\cU_{n,d}$. 

In \cite{EO}, we introduced an algebra $B'$ of pseudodifferential operators, 
quantizing the algebra $A'$ of Example \ref{ex:2}. Then 
$B'' = S^{\wt g}(B') = (B^{\prime \otimes \wt g})^{\SG_{\wt g}}$
is a quantization of $A''$. 

Then $Q(B)$ and $Q(B'')$ are quantizations of $R(\TT^* \cU_{n,d})$
and of $R(\KK^{[\wt g]})$, respectively.

\subsubsection{A quantization conjecture} \label{quant:conj}
 
It is natural to expect that the Beilinson-Drin\-feld adelic construction 
of twisted differential operators on the moduli stack $\on{Bun}_{n,d}$
yields a linear map $t : \cH^* \to \Diff(\cU_{n,d})$, whose image generates 
a commutative algebra. Once such a linear map is fixed, another such map 
can be obtained using any element $h\in \cH$: we set $(t+h)(\xi) = 
t(h) + \xi(h)1$. So we get an affine space of linear maps generating 
commutative algebras of differential operators; its underlying vector space
is $\cH$. 

On the other hand, we proposed in \cite{ER} the construction of a linear map 
$t'' : \cH^* \to Q(B'')$ with image generating a commutative algebra, 
starting from a linear map $\cH^* \to A'$. 

It is therefore natural to conjecture that for any linear map 
$$
t : \cH^* \to \Diff(\cU_{n,d}) = B
$$ 
provided by the Beilinson-Drinfeld
construction, one can find a linear map $\cH^* \to A'$ and an isomorphism 
$Q(B) \to Q(B'')$, such that the diagram 
$$
\begin{array}{ccccc}
Q(B) &         & \to &           & Q(B'') \\
    & \nwarrow &     & \nearrow   & \\
     &          & \cH^* & & 
\end{array}
$$
commutes and is a deformation of (\ref{comm:diag:fields}). 


\section{Explicit expressions in the rational case} \label{sect:1}

Theorem \ref{GNR} can be extended to the case of curves with 
marked points and parabolic bundles. In this section, we recall the 
explicit forms of the Hitchin and Beauville-Mukai systems in the case of 
$\PP^1$ with marked points and $n=2$ (Sections \ref{sect:1.1} and \ref{sect:1.2}). 
The rational version of the Hitchin system was introduced and studied
by Beauville in \cite{Acta}.  In Section \ref{sect:1.3}, we construct an 
explicit isomorphism between both systems and we explain in Section \ref{sect:1.4}
why it coincides with $\GNR$. 

\subsection{The rational Hitchin system} \label{sect:1.1}

\subsubsection{The fibration $\wt\cP_g \to \wt V_g$} 

If $n$ is an integer, we denote by $\cO(n)$ the line 
bundle of degree $n$ on $\PP^1$. We denote by $P_{\infty}$ the point at infinity of
$\PP^1$, and identify $\cO(n)$ with $\cO(nP_\infty)$.

Let $E = {\cO}\oplus{\cO(-1)}$. Then $E$ is a vector
bundle on $\PP^1$ of rank $2$ and degree $-1$. Let $g\geq 0$ 
be an integer. Let us set
$$
\Hom_{0}(E,E\otimes{\cO}(g+1)) = \{\phi\in
\Hom(E,E\otimes {\cO(g+1)})|\tr(\phi) = 0\} .
$$
\begin{defin}
Let us set
$$
{\wt{\cP}_g} = \Hom_{0}(E,E\otimes {\cO(g+1)})/\Aut(E), \; 
{\wt V_g} = H^0({\PP^1},{\cO(2g+2)}),  
$$
and let ${\wt H}_g : {\wt{\cP}_g}\to {\wt V_g}$ be the map 
induced by $\phi \mapsto {1\over2} \tr(\phi)^2.$
\end{defin}

Then we have the following description of ${\wt{\cP}_g}$ and ${\wt V_g}$. 
We set 
$$
{{\wt\cM}_g} = \{
M= \left(
\begin{array}{cc}
A& C\\
B& -A
\end{array}
\right) \in {\mathfrak sl}_2({\CC}[X]) \¦ 
\deg(A)\leq g, \deg(B)\leq g+1, \deg (C)\leq g+2
\} 
$$
($A,B,C$ are polynomials in one variable $X$) and 
$$
\Gamma = \{
 \left( \begin{array}{cc}
a& bX +c\\
0& a^{-1}
\end{array} \right) \¦ 
a\in {\CC}^*, b,c\in {\CC} \}
$$
($\Gamma$ is a subgroup of $\on{SL}_2(\CC[X])$). Then $\Gamma$ acts on 
${\wt\cM}_g$ by conjugation.

\begin{lemma}
${\wt{\cP}_g}$ identifies with the quotient ${{\wt\cM}_g}/{\Gamma},\, {\wt V_g}$
identifies with the space ${\CC}[X]_{\leq 2(g+1)}$ of polynomials of degree $\leq
2(g+1)$, and the map ${\wt H}_g$ identifies in its turn with the map induced by $M 
\mapsto A^2 + BC$.
\end{lemma}

If $P\in {\wt V_g}$, let us denote by $C_P$ the curve of $T^*(\PP^1)$, whose
intersection with $T^*({\PP^1}-\{P_{\infty}\})$ is defined by the equation 
$y^2 = P(x)$.

When $\deg(P) = 2(g+1)$, we have $\genus(C_P) = g$. According to \cite{Acta}, in this
case, the preimage ${\wt H}_{g}^{-1}(P)$ identifies with a Zariski open subset of
$\Jac^g(C_P)$.

\subsubsection{Poisson structure and integrable system on $\cP_g$}
 
Let $a_1,\ldots,a_{g+3}$ be distinct points of ${\PP^1}-\{P_{\infty}\}$.

\begin{lemma}
Set
$$
V_g = \{P\in {\wt V_g}|P(a_1) =\ldots=P(a_{g+3})= 0\}, \; 
{\cP}_g = {\wt H}_{g}^{-1}(V_g).
$$
The fibration ${\wt H}_{g} :{\wt {\cP}_g} \to {\wt V_g}$
restricts to a fibration $H_{g}^\Hitch :{\cP}_g \to V_g$. 
\end{lemma}

We now introduce a Poisson structure on ${\cP}_g$. For this, we will 
show that $\cP_g$ is isomorphic, as a variety, to a symplectic quotient
$\cN^g /\!/ \Gamma$, where $\cN$ is the nilpotent cone of $\SL_2(\CC)$, 
and transport on $\cP_g$ the Poisson structure of 
$\cN^g /\!/ \Gamma$. To construct the isomorphism $\cP_g \to \cN^g /\!/ \Gamma$, 
we will need the following explicit description of $\cP_g$. 

\begin{lemma}
Let us set
$$
{\cM}_g = \{M\in {\wt {\cM}_g}|\tr(M^2)(a_1) =\ldots=\tr(M^2)(a_{g+3})= 0\}.
$$
Then the action of $\Gamma$ preserves $\cM_g$, and $\cP_g$ identifies
with the quotient $\cM_g / \Gamma$.  
\end{lemma}

\medskip 

We now describe the symplectic quotient $\cN^{g+3}/\!/ \Gamma$. 
We will denote by ${\cN} = \{N\in {\mathfrak sl}_2(\CC)|\tr(N^2)=0\}$ the nilpotent
cone of ${\mathfrak sl}_2(\CC)$. The group $\Gamma$ acts on the product ${\cN}^{g+3}$
via
\begin{equation}
\label{gamma}
\gamma\cdot(A_1,\ldots,A_{g+3}):=
\big( \Ad(\gamma(a_1))(A_1),\ldots,\Ad(\gamma(a_{g+3}))(A_{g+3})\big),\,
\gamma \in \Gamma.
\end{equation}

Recall that the nilpotent cone $\cN$ (and hence ${\cN}^{g+3}$) has a
natural Poisson structure; we equip ${\cN}^{g+3}$ with the product structure. 
Then the action (\ref{gamma}) is hamiltonian, and its moment map is given by 
$$
\mu : {\cN}^{g+3}\to (\Lie \Gamma)^* = {\CC}^3, 
$$
$$
(A_1,\ldots,A_{g+3})\to
(\sum_{i=1}^{g+3} v_i,\sum_{i=1}^{g+3} u_i,\sum_{i=1}^{g+3}a_iu_i),
$$
where $(u_i,v_i,w_i)$ are the coordinates of $A_i \in {\mathfrak sl}_2(\CC)$ 
in the the Chevalley basis $e,f,h$, and we identify $\Lie \Gamma$ with $\CC^3$
using the basis $(h\otimes 1, e \otimes 1, e \otimes X)$. 

Then the symplectic quotient ${\cN}^{g+3} /\!/{\Gamma}$ is equal to the quotient 
$\mu^{-1}(0)/{\Gamma}$, which is equipped with a natural Poisson structure.

\medskip 

We now construct the isomorphism $\cP_g \to \cN^{g+3} /\!/ \Gamma$. 
For each $M\in \cM_g$, and each $i = 1,\ldots,g+3$, let us set 
$$
A_i(M) = \frac{M(a_i)}{\prod_{j\¦j\neq i}(a_i - a_j)}. 
$$
The family $(A_i(M))_{i = 1,\ldots,g+3}$ is uniquely determined by 
the identity 
$$
\frac{M(X)}{{\Pi}_{i=1}^{g+3}(X-a_i)} = \sum_{i=1}^{g+3}\frac{A_i(M)}{X-a_i}. 
$$
Then the assignment $M\mapsto (A_1(M),\ldots,A_{g+3}(M))$ defines a 
morphism
\begin{equation} \label{star}
i' : {\cM}_g \to {\cN}^{g+3}. 
\end{equation}
The image of this morphism is contained in $\mu^{-1}(0)$, and it actually 
induces an isomorphism $\cM_g \to \mu^{-1}(0)$. This morphism is 
equivariant with respect to the action of $\Gamma$.

We have therefore: 

\begin{prop} (see \cite{Acta})
$i'$ induces an isomorphism 
\begin{equation}\label{def:i4}
i: {\cP}_g = {\cM}_g /{\Gamma} \to \mu^{-1}(0)/{\Gamma}={\cN}^{g+3} /\!/{\Gamma}.
\end{equation}
of algebraic varieties, and we define a Poisson structure
on $\cP_g$ uniquely by the condition that $i$ is a Poisson morphism
(this structure is actually symplectic at a generic point). 
 
Moreover, the family of the functions on $\cP_g$ obtained
as $(H_{g}^\Hitch)^*(\ell)$ Poisson-commute ($\ell$ a function on $V_g$), 
so the fibration $H^\Hitch_{g} : {\cP}_g \to V_g$ is
Lagrangian. 
\end{prop}

This is in fact an algebraic completely integrable system. 
We will call this system the {\it rational Hitchin system}. 

\begin{remark}
We have 
$\dim ({\wt {\cP}_g})=(3g+6)-3=3g+3$, $\dim(\wt V_g)= 2g+3$,
and $\dim ({\cP}_g)=2g$, $\dim({V}_g)= g$.
\end{remark}
\medskip 

\subsubsection{Relation with the Hitchin systems} \label{rel:hitch} 

Let $N$ be arbitrary and let us set $\ul a = (a_1,\ldots,a_N)$. 
Hitchin's integrable system is 
defined on the cotangent bundle $T^*(\Bun_{({\PP^1},\ul a)}(n))$, where 
$\Bun_{({\PP^1},\ul a)}(n)$ is the moduli stack of vector bundles of
rank $n$ on $\PP^1$, with parabolic structures at $\underline{a}$.
$\Bun_{({\PP^1},\ul a)}(n)$ is the disjoint union of all 
$\Bun_{({\PP^1},\ul a)}(n,d)$, $d\in\ZZ$, where $\Bun_{({\PP^1},\ul a)}(n,d)$
is the moduli stack of parabolic structures on vector bundles of rank $n$
and of degree $d$.  The spaces $\Bun_{({\PP^1},\ul a)}(n,d)$ and 
$\Bun_{({\PP^1},\ul a)}(n,d+n)$ are isomorphic, so we assume $d\in [1-n,0]$. 

The open cell of the moduli space $\Bun_{({\PP^1},\ul a)}(n,d)$ corresponds 
to parabolic structures on the bundle ${\cal O}^{\oplus (n+d)} \oplus 
{\cal O}(-1)^{\oplus(-d)}$. The system studied here corresponds to $n=2,d=-1$. 

The reason we are studying the case $(n,d)= (2,-1)$ rather than $(n,d)= (2,0)$ 
is that in the latter case, the fibers of $(H^\Hitch)^{-1}(P)$ are isomorphic, 
when $P$ is generic, to the degree $\wt g-1$ Jacobian ${\Jac}^{\wt g-1}(C_P)$ 
(here $\wt g = \on{genus}(C_P)$, see \cite{Acta}); the generic line bundle 
of this Jacobian has no section. 
In the former case, the fibers are isomorphic to ${\Jac}^{\wt g}(C_P)$; generic 
bundles of this Jacobian have a unique (up to scalar multiplication) 
non-zero section and are uniquely determined by the set of zeroes of this section. 

In the general case, this situation corresponds to $d = 1-n$.

\subsubsection{Reformulation (the rational Hitchin system 
as a classical limit of the Gau\-din system.)} \label{sect:HHitch}

Let us reexpress the rational Hitchin system as a classical limit 
of the Gaudin magnetic model (\cite{Gaudin}).

For $i=1,\ldots,g+3$ let us define the functions $\ul H_i^\Hitch : 
{\cN}^{g+3} /\!/{\Gamma} \mapsto {\CC}$ by
$$
\ul H_i^\Hitch (\on{class\ of\ }(A_1,\ldots,A_{g+3})) 
= \sum_{j\neq i}\frac{\tr(A_i A_j)}{a_i -a_j},
$$
for any $(A_1,\ldots,A_{g+3}) \in \mu^{-1}(0)$. 

The map ${\underline H}^\Hitch= (\ul H^\Hitch_1,\ldots,\ul H^\Hitch_{g+3})$ 
maps ${\cN}^{g+3}/\!/{\Gamma}$ to the $g$-dimensional vector subspace 
${\underline{V}_g}$ of ${\CC}^{g+3}$, defined as the set of all 
$(\ul h_1,\ldots,\ul h_{g+3})$ obeying the equations
$$
\sum_{i=1}^{g+3} \ul h_i = \sum_{i=1}^{g+3}a_i \ul h_i = 
\sum_{i=1}^{g+3}a{_i}^2 \ul h_i = 0.
$$

\begin{lemma}
The fibrations ${\underline H}^\Hitch: {\cN}^{g+3}/\!/{\Gamma} \mapsto {\underline{V}_g}$
and $H^\Hitch : {\cP}_g \to V_g$ are isomorphic.
\end{lemma}

More precisely, we have ${\underline H}^\Hitch= i_{V}\circ H^\Hitch
\circ i^{-1}$, where $i_{V} : V_{g}\to {\underline{V}_g}$ 
is the linear isomorphism taking a polynomial
$P$ of $V_g$ to the vector $(P'(a_1),\ldots,P'(a_{g+3}))$; the inverse map takes
the vector $(\ul h_1,\ldots,\ul h_{g+3})$ to the polynomial 
$$
\prod_{i=1}^{g+3}(X - a_i)^2 \sum_{i =1}^{g+3}
\frac{\ul h_i}{X - a_i}.
$$

\subsection{The Beauville-Mukai system} \label{sect:1.2}

In \cite{Cortona}, Beauville introduced integrable systems associated with 
a symplectic surface and a linear system. This construction was later 
generalized to Poisson surfaces (see \cite{Bottacin}). This family of integrable
systems can be described as follows. Let $S$ be a Poisson
surface, let ${\CC}(S)$ be the function field of $S$ and let $g$ be an
integer. 

Then $({\CC}(S)^{\otimes g})^{\frak S_{g}}$ may be viewed as a ring of
functions over $S^{(g)}= S^g/{\frak S_{g}}$ and has a natural Poisson structure.

Let $f_1,\ldots,f_{g-1},f$ be elements of ${\CC}(S)$, such that the family
$(1,f_1,\ldots,f_{g-1},f)$ is free. For $\varphi \in {\CC}(S)$, let $\varphi^{(i)}$ be the
image of $\varphi$ in the $i$-th factor of $({\CC}(S)^{\otimes g})^{\frak S_{g}}$, so
$\varphi^{(i)} = 1^{\otimes i-1}\otimes\varphi\otimes 1^{\otimes n-i}.$

For $i= 1,\ldots, g$, let $M_i$ be the $i$-th minor of the matrix
$$
\left(
\begin{array}{ccccc}
1& f{_1}^{(1)}& \ldots& f{_{g-1}}^{(1)}& f^{(1)}\\
1& f{_1}^{(2)}& \ldots& f{_{g-1}}^{(2)}& f^{(2)}\\
\vdots & \vdots& \vdots &   \vdots &         \vdots \\
1& f{_1}^{(g)}& \ldots& f{_{g-1}}^{(g)}& f^{(g)}
\end{array}
\right)
$$
(obtained by removing the $i$-th column and taking the determinant) and set 
$$
H_i^{\BM} = (-1)^{g+1-i}\frac{M_i}{M_{g+1}}.
$$

Then $(H^\BM_1,\ldots,H^\BM_g)$ is a Poisson commuting family of functions on 
$S^{(g)}$ (in this generality, this statement is proved in \cite{ER}; it is a 
generalization of the "birational" part of the involutivity statements of 
\cite{Cortona,Bottacin}). 

Moreover, $(H^\BM_1,\ldots,H^\BM_g)$ defines a Lagrangian fibration
$$
H^\BM_S : S^{(g)} \supset U_S  \to {\CC}^g,
$$
where $U_S$ is a Zariski open subset in $S^{(g)}$.

We will be interested in the situation where $S = T^*({\PP^1})$. We consider a
collection $\underline{a} = (a_1,\ldots a_{g+3})$ of distinct points in $\CC$ 
and we identify $T^*({\PP^1} - \{P_\infty\})$ with $T^*({\CC})= {\CC}^2$. 
Then ${\CC(S)} =
{\CC}(X,Y)$ ($X,Y$ are algebraic independent variables, $X$ is the canonical 
coordinate on $\CC$ and $Y$ is the coordinate on the cotangent fibers). 
It is equipped with the Poisson bracket
$$
\{X,Y\} = -2 \prod_{i = 1}^{g+3}\left(X - a_{i}\right) .
$$

We set 
$$
f_1 = X, \; f_2 = X^2,\ldots, f_{g-1} = X^{g-1}, \;  
f = \frac{Y^2}{\prod_{i =1}^{g+3} (X - a_{i})}.
$$

The hamiltonians $H_k^\BM$ can be computed explicitly as follows: the variables are
$\ul X = (X_1,\ldots,X_g)$ and $\ul Y = (Y_1,\ldots,Y_g)$. For $k=1,\ldots g$, we set 
\begin{align*}
& H_{g+1-k}^\BM(\underline{X},\underline{Y}):=
\\& 
(-1)^k \sum_{\al=1}^{g}\frac{1}{\prod_{\beta \¦ \beta\neq \al}(X_\al - X_\beta)}
\frac{Y_{\al}^2}{\prod_{i =1}^{g+3} (X_\al - a_i)}
\sum_{  
\stackrel{(\beta_1,\ldots, \beta_k)| \beta_1 < \cdots < \beta_k,} 
{\beta_1,\ldots,\beta_k \in \{1,\ldots,\check\al,\ldots,g\}}  }
X_{\beta_1}\cdots X_{\beta_k}
\end{align*}
and the Poisson brackets in $\CC(S)^{\otimes g}$ are
$$
\{X_\al,Y_\beta\}= - 2 \prod_{i = 1}^{g+3}\left(X_\al - a_{i}\right)
\delta_{\al\beta}, \quad 
\{X_\al,X_\beta\}= \{Y_\al,Y_\beta\}=0.
$$

We denote by $U_{\underline{a}}$ the Zariski open subset of 
$(T^*(\PP^1))^{(g)}$ where each $H^\BM_i$ is regular, and
$$
H^\BM_{\underline{a}} = (H^\BM_1,\ldots,H^\BM_g) : 
(T^*(\PP^1))^{(g)}\supset U_{\underline{a}} \to {\CC}^g
$$
the corresponding morphism.

\subsection{A birational equivalence between the rational Hitchin and the 
Beauville-Mukai systems} \label{sect:1.3}

We will construct inverse birational morphisms
$$
\alpha : (\cN^{g+3}/\!/\Gamma) \supset U \to (T^*(\PP^1))^{(g)}
$$
and
$$
\beta :  (T^*(\PP^1))^{(g)} \supset V \to (\cN^{g+3}/\!/\Gamma) 
$$
and a linear isomorphism $\eta : {\underline {V}}_g \to {\CC}^g$, such that the
diagram
\begin{equation}
\label{diag}
\begin{array}{ccc}
(\cN^{g+3}/\!/\Gamma) \supset U & \stackrel{\ul H^\Hitch}{\longrightarrow} 
& \ul{V}_g \\
{\scriptstyle{\al}}\downarrow && \downarrow{\scriptstyle{\eta}} \\
(T^*(\PP^1))^{(g)} \supset U_{\underline{a}} & 
\stackrel{H_{\ul a}^\BM}{\longrightarrow} & \CC^g
\end{array}
\end{equation}
commutes (the morphism $\ul H^\Hitch$ was introduced in Section \ref{sect:HHitch}).

$(\alpha,\eta)$ and $(\beta,\eta^{-1})$ will therefore set up inverse 
birational equivalences between the rational Hitchin and the Beauville-Mukai
systems. 

\medskip 
{\underline {Construction of $\alpha$}:} let $U$ be a subset of $(\cN^{g+3}/\!/\Gamma)$ 
consisting
of all classes of $(g + 3)$-uples $(A_1,\ldots,A_{g+3})$ such that 
$\sum_{i \geq 1}^{g+3}(a_i)^2 u_i \neq 0$. For any such $(g + 3)$-uple, the
polynomial
$$
\prod_{i = 1}^{g+3}(X - a_i)\sum_{i = 1}^{g+3}\frac{u_i}{X - a_i}
$$
has exactly $g$ zeroes; we denote this set of zeroes by $\{x_1,\ldots,x_g\}$.

For $\alpha = 1,\ldots,g$, we set
\begin{equation} \label{thl}
y_{\alpha} = (\sum_{i = 1}^{g+3}\frac{v_i}{x_\al - a_i})
\prod_{i = 1}^{g+3}(x_\al - a_i).
\end{equation}
Then $\{(x_\al,y_\al), \al=1,\ldots,g\}$ belongs to $T^*(\PP^1)^{(g)}$.
We will show that it depends only on the 
$\Gamma$-orbit of $(A_1,\ldots,A_{g+3})$, and we set 
$$
\al(\on{class\ of\ }(A_1,\ldots,A_{g+3})) = \{(x_\al,y_\al), \al=1,\ldots,g\}.  
$$

\medskip 

{\underline {Construction of $\beta$}:} let us set $V = T^*(\PP^1 -
\{a_1,\ldots,a_{g+3},P_{\infty}\})^{(g)} = ((\CC - \{a_1,\ldots,a_{g+3}\})
\times \CC)^{(g)}$, then $\beta$ maps $\{(x_\al,y_\al), \al=1,\ldots,g\}$ 
to the class of $(A_1,\ldots,A_{g+3})\in {\cN}^{g+3}$, defined by
\begin{equation}
\label{Aif}
 u_i = \frac{\prod_{\al=1}^{g}(a_i - x_\al)}
{\prod_{j =1, j \neq i}^{g+3}(a_i - a_j)};
\end{equation}
\begin{equation}
\label{Aih}
 v_i = \sum_{\al=1}^{g} y_\al \cdot 
\frac{\prod_{\beta\neq \al}^{g}(a_i - x_\beta)}
{\prod_{\beta\neq \al}(x_\al - x_\beta)
\prod_{j=1, j \neq i}^{g+3}(a_i - a_j)};
\end{equation}
\begin{equation}
\label{Aie}
 w_i = -\frac{v_i^2}{u_i}.
\end{equation}
We will show that $(A_1,\ldots,A_{g+3})\in \mu^{-1}(0)$, so there is a 
well-defined map $\beta$ taking $\{(x_\al,y_\al),\al = 1,\ldots,g\}$ to 
the class of $(A_1,\ldots,A_{g+3})$. 

\medskip 

{\underline {Construction of $\eta$}:} $\ul V_g$ consists of all 
$(\ul h_1,\ldots,\ul h_{g+3}) \in \CC^{g+3}$, such that
$$
\sum_{i=1}^{g+3}\ul h_i = \sum_{i=1}^{g+3}a_i \ul h_i 
= \sum_{i=1}^{g+3}(a_i)^2 \ul h_i = 0. 
$$
$\eta$ maps such a $(g+3)$-uple to the vector 
$(h_{1}^{\underline{a}},\ldots,h_{g}^{\underline{a}})$ of $\CC^g$,
defined by the identity
$$
(\sum_{i =1}^{g+3}\frac{\ul h_{i}}{X - a_{i}})
\prod_{i =1}^{g+3}(X - a_{i}) = \sum_{i=1}^{g}
h_{i}^{\underline{a}}X^i.
$$

In other words, we have
$$
h_{i}^{\underline{a}} = (-1)^{g-i}
\sum_{j = 1}^{g+3} \ul h_j
\sum_{\stackrel{
\stackrel{(j_1,\ldots,j_{g+3-i}) \¦}{ j_1,\ldots,j_{g+3-i} \in 
\{1,\ldots,\check j,\ldots,g+3\}}} 
{j_1 < \ldots < j_{g+3-i}} }
a_{j_1}\ldots a_{j_{g+3-i}}
$$
and
$$
\ul h_{i} = \frac{\sum_{j=0}^{g-1} h_{j+1}^{\underline{a}}(a_{i})^j}
{\prod_{k =1,k \neq i}^{g+3}(a_{i}-a_{k})}.
$$

\begin{thm} \label{thm:birat:isom}
The maps $\alpha$ and $\beta$ are inverse to each other on Zariski-open
subsets, and when restricted to suitable such subsets, they induce
isomorphism between the fibrations 
$\ul H^\Hitch: (\cN^{g+3}/\!/\Gamma) \rightarrow \underline{V}_g$ and 
$H^\BM_{\underline{a}} : (T^*(\PP^1))^{(g)} \rightarrow {\CC}^g$.
In other words, the diagram 
(\ref{diag}) commutes. 
\end{thm}

\medskip 

We will study in the next sections the Poisson aspects of this statement
(Poisson isomorphism between Poisson rings.)

\medskip 

{\em Proof.} Let us first show that the map $\alpha$ is well-defined. The
action of the element
$$
\gamma =
 \left(
\begin{array}{cc}
a& bX +c\\
0& a^{-1}
\end{array}
\right)
$$
of $\Gamma$ on $(A_1,\ldots,A_{g+3})$ multiplies
$$
(\sum_{i =1}^{g+3}\frac{v_i}{X-a_{i}}) 
\prod_{i = 1}^{g+3}(X -a_{i})
$$
by a nonzero scalar, and it adds to
$$
( \sum_{i =1}^{g+3}\frac{v_i}{X- a_{i}})
\prod_{i = 1}^{g+3}(X -a_{i}),
$$
the polynomial 
$$
a^{-1}(bX+c)\sum_{i=1}^{g+3} \frac{u_i}{X - a_{i}})
(\prod_{i = 1}^{g+3}(X -a_{i}) , 
$$
which vanishes at each $x_\al$. So the set $\{(x_\al,y_\al),\al = 1,\ldots,g\}$ 
is not changed by the action of $\gamma$. This shows that $\alpha$ is well-defined.

Let us show that $\beta$ is well-defined. The relations (\ref{Aif}) and (\ref{Aih}) imply
\begin{equation}\label{first} 
\sum_{i =1}^{g+3}\frac{u_i}{X- a_{i}} =
(\sum_{i=1}^{g+3} a_i^2 u_i)  
\frac{\prod_{\alpha=1}^{g}(X - x_{\alpha})}{\prod_{i=1}^{g+3}(X - a_{i})},
\end{equation}
and
\begin{equation}\label{second} 
\sum_{i=1}^{g+3}\frac{v_i}{X - a_{i}} =
\frac{\prod_{\alpha=1}^{g}(X - x_{\alpha})}
{\prod_{i=1}^{g+3}(X - a_{i})}
\sum_{\alpha=1}^g \frac{\wt y_\al}{X - x_\al},
\end{equation}
where
\begin{equation} \label{third} 
\wt y_{\alpha} =
\frac{y_{\alpha}}{\prod_{\beta=1,\beta\neq\alpha}^g (x_\alpha - x_\beta )}.
\end{equation}
Looking at the behavior at infinity of identities (\ref{first}) and (\ref{second}), 
we get
$$
\sum_{i=1}^{g+3}u_i = \sum_{i=1}^{g+3}a_i u_i = 0
$$
and $\sum_{i=1}^{g+3}v_i = 0$, so $(A_1,\ldots,A_{g+3})\in \mu^{-1}(0)$ and
$\beta$ is well-defined.

The maps $\alpha$ and $\beta$ are inverse to each other, because the equations
(\ref{first}), (\ref{second}), (\ref{third}) and (\ref{Aie}) uniquely define
$\{(x_{\alpha},y_{\alpha}), \alpha = 1,\ldots,g\}$ in terms of the class
$(A_1,\ldots,A_{g+3})$ and vice versa.

\medskip 

To complete the proof of Theorem \ref{thm:birat:isom}, 
let us prove that the diagram (\ref{diag}) commutes. Assume that
$\{(x_\al,y_\al),\al = 1,\ldots,g\}$ and the class of $(A_1,\ldots,A_{g+3})$ 
are related by (\ref{first}), (\ref{second}), (\ref{third}) and 
(\ref{Aie}). Then
$$
{1\over 2}\tr\left( \sum_{i=1}^{g+3}\frac{A_i}{X - a_i}\right)^2 =
\sum_{i=1}^{g+3}\frac{\underline{H}_i^{\Hitch}(A_1,\ldots,A_{g+3})}{X - a_i}.
$$
Set $A(X) = \sum_{i=1}^{g+3}\frac{A_i}{X - a_i}$, then the relations 
(\ref{first}) and (\ref{second}) imply that $A(x_\al)$ has the form 
$$
A(x_\al) = \left(
\begin{array}{cc}
y_\al & *\\
0& -y_\al
\end{array}
\right)\prod_{i=1}^{g+3}(x_\al - a_i)^{-1}
$$
so
$$
{1\over 2} \tr\left(  \sum_{i=1}^{g+3}\frac{A_i}{X - a_i}\right)^2 |_{X=x_\al} 
= \left(\frac{y_\al}{\prod_{i}(x_\al -a_i)}\right)^2.
$$
The behavior of $\sum_{i=1}^{g+3}\frac{A_i}{x - a_i}$ when $x\to \infty$ shows that
$$
\prod_{i=1}^{g+3}(X -a_i)\cdot {1\over 2} \tr\left(  
\sum_{i=1}^{g+3}\frac{A_i}{X - a_i}\right)^2
$$
is a polynomial of degree $\leq g-1$, so it is equal to
$$
\big( \prod_{\al=1}^g(X - x_\al) \big) 
\sum_{\al = 1}^g \frac{y_\al ^2}{X - x_\al}
\big(\prod_{\beta=1,\beta\neq\al}^{g}(x_\al -
x_\beta)\big)^{-1} \big(\prod_{i=1}^{g+3}(x_\al - a_i)\big)^{-1}.
$$
So we get
\begin{align*}
& \sum_{i=1}^{g+3}\frac{\ul H_i^\Hitch (A_1,\ldots,A_{g+3})}{X - a_i}
= \\ & 
\frac{\prod_{\al=1}^g (X - x_\al)}{\prod_{i=1}^{g+3} (X - a_i) }
\sum_{\al = 1}^g
\frac{y_\al ^2}{X - x_\al}\big(\prod_{\beta=1,\beta\neq\al}^{g}(x_\al -
x_\beta)\big)^{-1}\big(\prod_{i=1}^{g+3}(x_\al - a_i)\big)^{-1}.
\end{align*}
and it follows that
\begin{align*}
& \ul H_i^\Hitch (A_1,\ldots,A_{g+3})=
\\& 
\frac{\prod_{\al=1}^g(a_i - x_\al)}{\prod_{j=1,j\neq i}^{g+3}(a_i -
a_j)}\sum_{\al=1}^g\frac{y_\al ^2}{a_i -
x_\al}\big(\prod_{\beta=1,\beta\neq\al}^{g}(x_\al -
x_\beta)\big)^{-1}\big(\prod_{i=1}^{g+3}(x_\al - a_i)\big)^{-1}, 
\end{align*}
therefore 
$$
\ul H_i^\Hitch (A_1,\ldots,A_{g+3})=
\frac{1}{\prod_{j|j\neq i}(a_i - a_j)}\sum_{k=0}^{g-1}a_i ^k
H^\BM_{k+1}(\underline{x},\underline{y}),
$$
which shows the commutativity of (\ref{diag}). This ends the proof of 
Theorem \ref{thm:birat:isom}. 
\hfill \qed\medskip

\subsection{Relation with spectral transformation and $\GNR$} 
\label{sect:1.4}

Let us show that the map $\al$ coincides with the map $\GNR$. For this, 
we compute $\GNR$. 

Recall that ${\wt{\cP_g}} =
\Hom_0(E,E\otimes\cO(g+1))/\Aut(E)$ identifies with
${\wt{\cM}_g}/{\Gamma}$, where
$$
{\wt{\cM}_g}= \{\left(
\begin{array}{cc}
A& C\\
B& -A
\end{array}
\right), A,B\in{\CC}[X], \deg(B)\leq g, \deg(A)\leq g+1, \deg(C)\leq g+2 \}
$$
and
$$
{\Gamma}= \{\left(
\begin{array}{cc}
a& bX + c\\
0& a^{-1}
\end{array}
\right), a\in{\CC}^*, b,c\in{\CC}\}.
$$
In the definition of $\wt{\cP}_g$, the bundle $E$ may be replaced 
by $E \otimes \cO(k)$, where $k$ is any integer.

If $M\in \wt{\cM}_g$, the spectral curve $C(M)$ is defined by $y^2 = (A^2 +
BC)(x)$ (equation in $T^*(\PP^1)$). Assume that $\deg(A^2 + BC)= 2(g+1)$,
and write $(A^2 +BC)(X) = a_0^2 X^{2(g+1)} + \on{terms\ of\ degree\ }<2(g+1)$. 
Then the curve $C(M)$
has two branches at infinity: a local parameter at infinity is $x^{-1}$, and the
branches are defined by $y = \pm a_0 x^{g+1} + o(x^{g+1})$. We have a natural
projection $\pi : C(M)\mapsto \PP^1$.

If $M\in {\widetilde{\cM}_g}$, its {\it spectral bundle} is the line bundle over
$C(M)$, defined as the subbundle of $\pi^*(E \otimes \cO(k))$
whose fiber over $(x,y)\in C(M)$ is 
$$
{\cL}_M(x,y) = \Ker\left(M(x) -y\right)\subset \pi^*(E \otimes \cO(k))(x,y)
$$
(see \cite{Hitchin}).

When $k = g+1$, then $H^0(C(M),{\cL}_M)$ is one-dimensional, spanned by 
$$
\sigma(x,y)
=\left(\begin{array}{c}y + A(x)\\ B(x)\end{array}\right);  
$$
$\sigma(x,y)$ is regular at all points of $C(M)$ at finite distance, 
and at the branches at infinity, it has the expansion 
$$\left(\begin{array}{c}(\al\pm a_0)x^{g+1} + o(x^{g+1})\\
\beta x^g + o(x^g)\end{array}\right), 
$$
so it is a regular section of $E\otimes \cO(g+1)$
(here we set $A(x) = \al x^{g+1} + O(x^g)$, $B(x) = \beta x^g + o(x^g)$). 
Moreover, $\cL_M\in \Jac^{g}(C(M))$ belongs to the image of the injection 
$\AJ_g : C(M)^{(g)} \to \Jac^{g}(C(M))$. The preimage of $\cL_M$ in 
$C(M)^{(g)}$ is the collection of zeroes of $\sigma(x,y)$, which is the set 
of all $(x_i,y_i)$ defined as follows:

(1) $x_1,\ldots,x_g$ is the set of zeroes of $B(X)$;

(2) for $i = 1,\ldots,g$, $y_i = A(x_i)$. 

Then the map $\GNR$ associates to $M$ the collection 
$\{(x_i,y_i),i=1,\ldots,g\}$. So $\GNR = \al$. 


\section{Ringed Lie algebras and their reductions}
\label{sect:ringed}  

Our next task is to express explicitly the isomorphism between the Poisson 
fields of $\cN^{g+3}/\!/\Gamma$ and $\KK_{\PP^1}^{[g]}$. It turns out that these
Poisson fields are the fraction fields of graded subrings $S^\bullet_{\cX^0}(\cX^1)$
and $S^\bullet_{\cU^0}(\cU^1)$, which correspond to "ringed Lie algebras"
$\cX$ and $\cU$. In this section, we explain the formalism of ringed Lie 
algebras and the consturction of $\cX$ and $\cU$. 

\subsection{Ringed Lie algebras}

Let $R$ be a commutative ring over $\CC$, and let $V$ be a vector space over 
$\CC$. Then $R \otimes S^\bullet(V)$ is a graded algebra, where elements of
$R$ have degree $0$ and elements of $V$ have degree $1$. 

Set $\cG^0 = R$, $\cG^1 = V \otimes R$. Set $\cG = \cG^0 \oplus \cG^1$, and 
view $\cG$ as a graded complex vector space. It is also a graded 
$R$-module. 
Then if 
$$
\{,\} : \wedge^2(\cG) \to \cG 
$$
is a linear map of degree $-1$, then $\{r,f\} \in \cG^0 = R$ if $r\in R$
and $f\in \cG$. 

\begin{lemma}
There is a bijective correspondence between: 

(1) Poisson structures  on $R \otimes S^\bullet(V)$ of degree $-1$, 
and

(2) Lie brackets $\{,\} : \wedge^2( \cG) \to \cG$ of degree $-1$, 
satisfying the relations 
$$
\{rf,f'\} = r \{f,f'\} + \{r,f'\}f
$$ 
for any $r\in R$ and $f,f'\in \cG$. 
\end{lemma}

Such Lie brackets yield, and are uniquely determined by a pair of maps 
$\wedge^2(V)$ $\to (V \otimes R)$ and $V \otimes R \to R$. We say
that $(\cG, \{,\})$ is a 
{\it ringed Lie algebra}. 

\subsection{Reductions of ringed Lie algebras}

\subsubsection{Reduction of type I} \label{sect:red:I}

Let $R$ be a ring and let $\big(R \oplus (V \otimes R),\{,\}\big)$
be a ringed Lie algebra. Let $I$ be a prime ideal of $R$. 
Let us set 
$$
N = \{f\in V \otimes R \¦ \{f,I\} \subset I\}
$$
and $\cV = N / (V\otimes I)$. Then $\cV$ is a vector space over 
$R/I$. 

\begin{lemma}
$(R/I) \oplus \cV$ has a ringed Lie algebra stucture.  
\end{lemma}

{\em Proof.} We have a sequence of inclusions
$$
I \oplus (V \otimes I) \hookrightarrow R \oplus N 
\hookrightarrow  R \oplus (V\otimes R).  
$$
The first map is a Lie ideal inclusion, and the second map 
is an inclusion of Lie algebras. So $(R/I) \oplus \cV$ is a Lie 
algebra. One checks that it satisfies the Leibniz rules with respect to 
the structure of module over $R/I$. 
\hfill \qed \medskip 

\subsubsection{Reduction of type II}

Assume that $R$ is a field and let $\big(R \oplus (V \otimes R),\{,\}\big)$
be a ringed Lie algebra. Let $\Phi \in R\otimes V$ be a fixed element. 
Set $R^\Phi = \{r\in R \¦ \{\Phi,r\} = 0\}$, so $R^\Phi$ is a subfield of 
$R$, and $(R\otimes V)^\Phi = \{ f\in R\otimes V \¦ \{\Phi,f\} = 0\}$. 
Then $\{,\}$ restricts to $R^\Phi \oplus (R \otimes V)^\Phi$, and 
$\big(R^\Phi \oplus (R \otimes V)^\Phi, \{,\} \big)$ is a ringed Lie algebra.
Let us denote by $R^\Phi \Phi$ the $R^\Phi$-vector subspace of $(R \otimes V)^\Phi$
generated by $\Phi$. Then $\{0\} \oplus R^\Phi \Phi$ is a Lie ideal of 
$R^\Phi \oplus (R \otimes V)^\Phi$. Moreover: 

\begin{lemma}
$R^\Phi \oplus \big( (R\otimes V)^\Phi/ R^\Phi \Phi \big)$ 
is a ringed Lie algebra. 
\end{lemma}

\subsection{Other constructions of ringed Lie algebras}

\

\medskip \noindent 
(1) $R$ is a field, $\G$ is a $\CC$-Lie algebra, and we have a Lie 
algebra morphism  $\G \to \on{Der}(R)$. Then 
\begin{equation} \label{basic:constr}
R \oplus (\G\otimes R)
\end{equation} 
is a ringed Lie algebra. 

\medskip \noindent
(2) $R$ is a field, $R \oplus (V\otimes R)$ is a ringed Lie 
algebra, and $\Gamma$ is a finite group of automorphisms of 
$R \oplus (V\otimes R)$. Then $R^\Gamma$ is a field and 
$R^\Gamma \oplus (V\otimes R)^\Gamma$ is a ringed Lie 
algebra. 

\subsection{Examples}

\subsubsection{The ringed Lie algebra $\cX$}

We set $R = \CC(x_\al,\al = 1,\ldots,g)$; $\G$ is the abelian Lie algebra 
$\oplus_{\al = 1}^g \CC y_\al$, and the morphism $\G \to \on{Der}(R)$
takes $y_\al$ to $ - \prod_{i = 1}^{g+3} (x_\al - a_i) {\pa\over{\pa x_\al}}$. 
We therefore obtain a ringed Lie algebra 
$$
\cX' = \CC(x_\al,\al = 1,\ldots,g) \oplus 
\big( (\oplus_{\al = 1}^g \CC y_\al) \otimes 
\CC(x_\al,\al = 1,\ldots,g) \big) . 
$$

Then $\Gamma = \SG_g$ acts on this ringed Lie algebra by permuting the 
pairs $(x_\al,y_\al)$. We set $\cX = (\cX')^\Gamma$. Then we have 
$$
\cX = \CC(x_\al,\al = 1,\ldots,g)^{\SG_{g}} \oplus 
\CC(x_\al,\al = 1,\ldots,g)^{\SG_{g-1}} ; 
$$
in the last vector space, $\SG_{g-1}$ acts by permuting the last 
variables $x_2,\ldots,x_{g}$. This is a vector space over 
$\CC(x_\al,\al = 1,\ldots,g)^{\SG_g}$
($g$-dimensional by Galois theory). An element 
$f(x_1,\ldots,x_g)$ of $\CC(x_\al,\al = 1,\ldots,g)^{\SG_{g-1}}$
corresponds to the element 
$$
Y[f] = \sum_{\al = 1}^g f(x_\al,x_1,\ldots,\check x_\al,\ldots,x_g) y_\al
$$
of $(\cX')^\Gamma$.

\subsubsection{Ringed Lie algebras associated with vector spaces}
\label{sect:vector}

Let $V$ be a finite-dimen\-sio\-nal vector space, and let $W\subset V$ be a subspace. 
Let $\CC(V)_W$ be the localization of the symmetric algebra $S(V)$
with respect to $V \setminus W$. Let $I \subset \CC(V)_W$ be the ideal 
$(W)$. Then $\CC(V)_W / I = \CC(V/W)$. 

Let $\cG(V,W)$ be the ringed Lie algebra
$$
\cG(V,W) = \CC(V)_W \oplus \big( \CC(V)_W \otimes V^*\big),   
$$
where the brackets are $\{\xi,v\}= \langle \xi,v\rangle$, 
$\{\xi,\xi'\} = 0$ for $\xi,\xi'\in V^*$, $v\in V$. 

Then the reduction of $\cG(V,W)$ with respect to $(W)$ is 
\begin{equation} \label{G:V:W}
\cG(V/W,0) = \CC(V/W) \oplus \big( \CC(V/W) \otimes (V/W)^*\big).   
\end{equation}
We have 
$$
N = \big( \CC(V)_W \otimes W^\perp\big) + (W) \otimes V^*.  
$$

The reduction of $\cG(V,0)$ with respect to $\Phi = \sum_i x_i \otimes \xi^i = $
the canonical element of $V \otimes V^*$ is 
\begin{equation} \label{G:V:0}
\CC(V)^0 \oplus \big( (\CC(V)^0 \otimes V^*) / \CC(V)^{-1} \Phi \big) , 
\end{equation} 
where $\CC(V)^0$ (resp., $\CC(V)^{-1}$) is the field (resp., space) of 
all rational functions on $V$, homogeneous of degree $0$ (resp., $-1$). 
Note that (\ref{G:V:0}) is isomorphic to the algebra $\cG(H,0)$ , where $H$
is any hyperplane of $V$, not containing the origin. 

Applying this reduction to $\cG(V/W,0)$, we obtain the ringed 
Lie algebra 
\begin{equation} \label{reduced:GLA}
\cU = \cU^0 \oplus \cU^1, 
\end{equation}
where $\cU^0 = \CC(V/W)^0$ and 
\begin{align*}
\cU^1 & = N^\Phi / \big( (I \otimes V^*)^\Phi + \CC(V/W)^{-1} \Phi \big) 
\\ & = \big( \CC(V/W)^0 \otimes (V/W)^* \big) / \big( \CC(V/W)^{-1} \Phi \big).  
\end{align*}

\subsubsection{The ringed Lie algebra $\cU$} \label{3.4.2}

We apply the constructions of Section \ref{sect:vector} to 
$$
V = \oplus_{i=1}^{g+3} \CC u_i, \; 
W = \on{Span}(\ell_0,\ell_1), 
$$
where $\ell_0 = \sum_i u_i$ and $\ell_1 = \sum_i a_i u_i$. 
We denote by $(u_i^*)_{i = 1,\ldots,g+3}$ the basis of $V^*$ dual to 
$(u_i)_{i = 1,\ldots,g+3}$, and we set $v_i = -2u_i u_i^*$. 

The relation of this notation with that of Section \ref{sect:1} is 
$$
u_i = (A_i)_f, \quad v_i = (A_i)_h.  
$$

Then $\cU^0 = K$, where $K$ is the field 
$$
\CC(\ol u_1,\ldots,\ol u_{g+3} \¦ \sum_i \ol u_i = 
\sum_i a_i \ol u_i = 0, \sum_i a_i^2 \ol u_i = 1), 
$$
and 
$$
\cU^1 = N^\Phi / \big( (I \otimes V^*)^\Phi + \CC(V/W)^{-1} \Phi \big) .  
$$

Here $\ol u_i$ is the class of $u_i / (\sum_i a_i^2 u_i)$. 
We use the following notation: if $x_1,\ldots,x_N$ are indeterminates, 
and $\ell_1,\ldots,\ell_k$ is a collection of independent linear combinations
of $x_1,\ldots,x_N$, then the field $\CC(x_1,\ldots,x_N\¦\ell_1 = a_1,\ldots,
\ell_k = a_k)$ is the function field of the affine subspace defined by
the equations $\ell_1 = a_1,\ldots,\ell_k = a_k$. 

On the other hand, $\cU^1$ is a $g$-dimensional $K$-vector space. 
We have 
$$
\cU^1 \simeq \{(\la_1,\ldots,\la_{g+3}) \in K^{g+3} \¦ \sum_i \la_i 
= \sum_i a_i \la_i = 0\} / K (\ol u_1,\ldots,\ol u_{g+3}). 
$$
The inverse isomorphism takes the class of $(\la_1,\ldots,\la_{g+3})$
to the class of 

\noindent $\sum_i \la_i (v_i / \ol u_i)$. 


\section{An isomorphism of ringed Lie algebras}
\label{sect:iso:ringed}  

We have $\CC(V)_W = \CC(u_1,\ldots,u_{g+3})_{\on{Span}(\ell_0,\ell_1)}$, 
where $\ell_0 = \sum_i u_i$, $\ell_1 = \sum_i a_i u_i$. Define 
$P_u(X)\in \CC(V)_W[X]$ as 
$$
P_u(X) = \prod_{i=1}^{g+3} (X-a_i) \cdot 
\sum_{i=1}^{g+3} {{u_i}\over{X-a_i}}. 
$$

\begin{thm} \label{thm:iso:class}
1) For any $\phi\in\CC[X]$, the element 
$$
W[\phi] = \on{Res}_{t = \infty} \Big( \sum_i {v_i \over{t-a_i}}
\big( {{P'_{u}}\over{P_{u}}}(t) 
- {{P'_{u}}\over{P_{u}}}(a_i)\big) 
\phi(t) dt \Big)  
$$
belongs to $N^\Phi$ and therefore defines an element $W[\phi] \in \cU^1$
(recall that $v_i = -2 u_i u_i^*$). 

2) There is a unique ringed Lie algebra morphism 
$\La : \cX \to \cU$, 
whose degree zero component $\cX^0 \to \cU^0$ is defined by 
\begin{equation} \label{im:sigma}
\sum_{\al_1 < \ldots < \al_k} x_{\al_1} \cdots x_{\al_k}
\mapsto \sum_{i=1}^{g+3} \ol u_i \sum_{
\stackrel{j_1 < \cdots < j_{k+2}}{j_1,\ldots,j_{k+2} \neq i}}
a_{j_1} \cdots a_{j_k} 
\end{equation} 
and whose degree $1$ component $\cX^1 \to \cU^1$ is defined by  
\begin{equation} \label{im:Y}
Y[\phi] = \sum_\al (\phi/\Pi)(x_\al) y_\al \mapsto 
W[\phi] = \on{Res}_{t = \infty} \Big( \sum_i {v_i \over{t-a_i}}
\big( {{P'_{u}}\over{P_{u}}}(t) 
- {{P'_{u}}\over{P_{u}}}(a_i)\big) 
\phi(t) dt \Big)  
\end{equation} 
where $\phi$ is any element of $\CC[X]$.  
Here $\Pi(X) = \prod_{i=1}^{g+3} (X-a_i)$.
\end{thm}

{\em Proof.}
Let us prove 1). Let $n$ be an integer $\geq 0$. 
$W[X^n]$ is the image of the coefficient of $X^{-n-1}$ in 
\begin{equation} \label{expr'}
\sum_{i=1}^{g+3} {{v_i}\over{X-a_i}} \big( 
{{P'_u}\over{P_u}}(X) - {{P'_u}\over{P_u}}(a_i)\big)
\in X^{-1} \cU^1[[X^{-1}]].   
\end{equation}
Let us compute the Poisson bracket of (\ref{expr'}) with $\ell_0$ and 
$\ell_1$. We find 
\begin{align*} 
\{\ell_0, (\ref{expr'}) \}& = 2 \sum_{i=1}^{g+3} {{u_i}\over{X-a_i}}
\big( {{P'_u}\over{P_u}}(X) - {{P'_u}\over{P_u}}(a_i)\big) 
\nonumber \\& = 
2 {{P'_u(X)}\over{\prod_i (X-a_i)}} - 2 \sum_i 
{{u_i}\over{X-a_i}} {{P'_u}\over{P_u}}(a_i) = 0;    
\end{align*}
the last equality is 
$$
{{P'_u(X)}\over{\prod_i (X-a_i)}} = \sum_i {1\over{X-a_i}}
u_i {{P'_u}\over{P_u}}(a_i), 
$$
which follow from simple poles decomposition of the left side and 
the identity $P_u(a_i) = u_i \prod_{j\¦j\neq i}(a_i - a_j)$.  
Therefore 
$$
\{\ell_0,W[X^n]\} = 0. 
$$
Next, 
\begin{align}
\nonumber 
& \{\ell_1,(\ref{expr'})\} = 2 \sum_{i=1}^{g+3} 
{{a_iu_i}\over{X-a_i}} \big( {{P'_u}\over{P_u}}(X) - 
{{P'_u}\over{P_u}}(a_i)\big) 
\\ & = \nonumber 
- 2 \big( \sum_i u_i\big) {{P'_u}\over{P_u}}(X) + 2 \sum_i u_i  
{{P'_u}\over{P_u}}(a_i)
+ 2 X \big( {{P'_u(X)}\over{\prod_i (X-a_i)}} 
- \sum_i {{u_i}\over{X-a_i}} {{P'_u}\over{P_u}}(a_i)\big) 
\\ & = \label{rhs:expr:2} 
- 2 \big( \sum_i u_i\big) {{P'_u}\over{P_u}}(X) + 2 \sum_i u_i  
{{P'_u}\over{P_u}}(a_i) .  
\end{align}
Now the first term of (\ref{rhs:expr:2}) belongs to $X^{-1} I[[X^{-1}]]$
($I$ is the ideal $(W)\subset \CC(V)_W$) and the second term belongs to $\CC(V)_W$. 
Therefore 
$$
\sum_i u_i  {{P'_u}\over{P_u}}(a_i) = 0,   
$$
and $W[X^n]\in N$. Clearly, $\{\Phi,W[X^n]\} = 0$, so $W[X^n] \in N^\Phi$.   
This proves 1). 
 
\medskip 

2) To prove that $\La$ is a linear isomorphism, we construct its inverse
$\La^{-1}$. 
It is given by the formulas 
$$
(\La^{-1})^0 : \ol u_i \mapsto {{\prod_{\al = 1}^g (a_i - x_\al)}\over 
{\prod_{j\¦j\neq i} (a_i - a_j)}}
$$
for the degree zero part, and 
$$
(\La^{-1})^1 : 
\on{class\ of\ }\sum_i \la_i (v_i /\ol u_i)
\mapsto 
- \sum_\al  {{y_\al / \Pi(x_\al)}
\over{\prod_{\beta \¦ \beta \neq \al} (x_\al - x_\beta)}}
\big( \sum_i {{(\La^{-1})^0(\la_i) }\over{x_\al - a_i}} \big) 
$$
for the degree $1$ part, where $(\la_1,\ldots,\la_{g+3})\in K^{g+3}$
is such that $\sum_i \la_i = \sum_i a_i \la_i = 0$. 

\medskip 

Let us now prove that $\La$ is a morphism of ringed Lie algebras. 
We introduce the generating series 
$$
X(t) = \sum_\al {1\over{t-x_\al}} \in \cX^0[[t^{-1}]], \; 
Y(t) = \sum_{\al} {{y_\al / \Pi(x_\al)}\over{t-x_\al}} \in 
\cX^1[[t^{-1}]] , 
$$
and 
$$
U(t) = {{P'_u(t)}\over{P_u(t)}} \in \cU^0[[t^{-1}]], \; 
W(t) = \sum_i {{v_i}\over{t-a_i}} \big( {{P'_u}\over{P_u}}(t)
- {{P'_u}\over{P_u}}(a_i)\big) \in \cU^1[[t^{-1}]].  
$$

Then the relations $\{x_\al,y_\beta\} = -2\Pi(x_\al)\delta_{\al,\beta}$
and $\{y_\al,y_\beta\} = 0$ imply 
\begin{equation} \label{PB:1}
\{Y(t),X(s)\} = 2 \pa_s \Big( {{X(s) - X(t)}\over{s-t}} \Big)  
\end{equation}
and 
\begin{equation} \label{PB:2}
\{Y(t),Y(s)\} = {2\over{s-t}} \Big( \pa_s Y(s) + \pa_t Y(t) 
- {2\over{s-t}} \big( Y(s) - Y(t) \big) \Big).  
\end{equation}
These relations generate the brackets of $\cX$. 

On the other hand, we have $\La(X(t)) = U(t)$ and $\La(Y(t)) = W(t)$. 
So it suffices to prove that relations (\ref{PB:1}) and (\ref{PB:2})
are satisfied, when $X(t),Y(t)$ are replaced by $U(t),W(t)$. 
We prove this in the next two lemmas. 
 
\begin{lemma} \label{lemma:identity}
\begin{equation} \label{PB:1bis}
\{W(t), U(s)\} = 
\pa_s \big( {2\over{s-t}} \big( U(t) - U(s)\big) \big).
\end{equation}
\end{lemma}

{\em Proof.} We will prove 
$$
\{W(t), \log(P_u(s))\} = 
{2\over{s-t}} \big( U(t) - U(s)\big) .
$$
Then the derivative of this identity with respect to $s$ is 
(\ref{PB:1bis}). 

Since 
$$
\{v_i, \log(P_u(s))\} = -2 {{u_i}\over{s-a_i}} {{\Pi(s)}\over{P_u(s)}}, 
$$
we find 
\begin{align} \label{interm} 
 & \nonumber \{\sum_i  {{v_i}\over{t-a_i}}
\big( {{P'_u}\over{P_u}}(t) - {{P'_u}\over{P_u}}(a_i)\big) 
, \log(P_u(s))\} 
\\& \nonumber = 
-2 \sum_i {{\Pi(s)}\over{P_u(s)}} {{u_i}\over{(t-a_i)(s-a_i)}}
\big( {{P'_u}\over{P_u}}(t) - {{P'_u}\over{P_u}}(a_i) \big) 
\\& \nonumber = 
- {2\over{s-t}} \sum_i {{\Pi(s)}\over{P_u(s)}} 
\big( {{u_i}\over{t-a_i}} - {{u_i}\over{s-a_i}} \big) 
\big( {{P'_u}\over{P_u}}(t) - {{P'_u}\over{P_u}}(a_i) \big) 
\\& \nonumber 
= - {2\over{s-t}} {{\Pi(s)}\over{P_u(s)}} 
{{P'_u(t)}\over{P_u(t)}} 
\big( {{P_u(t)}\over{\Pi(t)}} - {{P_u(s)}\over{\Pi(s)}} \big) 
+ {2\over{s-t}} {{\Pi(s)}\over{P_u(s)}} 
\sum_i \big( {{u_i}\over{t-a_i}} - {{u_i}\over{s-a_i}}\big) 
{{P'_u}\over{P_u}}(a_i)
\\&  =
{2\over{s-t}} {{P'_u(t)}\over{P_u(t)}} 
- {2\over{s-t}} {{\Pi(s)}\over{P_u(s)}} 
\Big( {{P'_u(t)}\over{\Pi(t)}} - \sum_i \big( {{u_i}\over{t-a_i}}
- {{u_i}\over{s-a_i}} \big) {{P'_u}\over{P_u}}(a_i)
\Big) . 
\end{align}
On the other hand, the polynomials (in $t$) $P'_u(t)$ and 
$$
\Pi(t) \sum_i {{u_i}\over{t-a_i}} {{P'_u(a_i)}\over{P_u(a_i)}} 
$$ 
both have degree $\leq g+2$, and since $P_u(a_i) = u_i \Pi'(a_i)$, 
their values for $t = a_i$ coincide for each $i = 1,\ldots,g+3$. 
Therefore 
$$
P'_u(t) = \Pi(t) \sum_i {{u_i}\over{t-a_i}} {{P'_u(a_i)}\over{P_u(a_i)}} ; 
$$
it follows that 
$$
{{P'_u(t)}\over{\Pi(t)}} - \sum_i \big( {{u_i}\over{t-a_i}}
- {{u_i}\over{s-a_i}} \big) {{P'_u}\over{P_u}}(a_i)
= {{P'_u(s)}\over{\Pi(s)}};   
$$
plugging this identity into (\ref{interm}), we get 
$$
\{\sum_i  {{v_i}\over{t-a_i}}
\big( {{P'_u}\over{P_u}}(t) - {{P'_u}\over{P_u}}(a_i)\big) 
, \log(P_u(s))\} = {2\over{s-t}} \big( {{P'_u(t)}\over{P_u(t)}}
- {{P'_u(s)}\over{P_u(s)}} \big).  
$$
This ends the proof of Lemma \ref{lemma:identity}. 
\hfill \qed \medskip 

We now prove 
\begin{lemma} \label{lemma:PB2}
We have 
\begin{equation} \label{PB:2bis}
\{W(t),W(s)\} = {2\over{s-t}} \Big( \pa_s W(s) + \pa_t W(t) 
- {2\over{s-t}} \big( W(s) - W(t) \big) \Big).  
\end{equation}
\end{lemma}

{\em Proof.} We have 
$$
\{v(s),u(t)\} = 2{{u(s) - u(t)}\over{s-t}}, \; \on{so} \; 
\{v(s),{{P'_u}\over{P_u}}(t)\} = \pa_t \Big( {2\over{s-t}} \big( 
{{\Pi(t)}\over{\Pi(s)}} {{P_u(s)}\over{P_u(t)}} - 1 \big) \Big) ,  
$$
therefore 
$$
\{v(s) {{P'_u}\over{P_u}}(s), v(t) {{P'_u}\over{P_u}}(t)\}
= v(t){{P'_u}\over{P_u}}(s) \pa_t \Big( {2\over{s-t}} \big( 
{{\Pi(t)}\over{\Pi(s)}} {{P_u(s)}\over{P_u(t)}} - 1 \big) \Big) - 
(s\leftrightarrow t).   
$$
Taking the singular part of both sides at $s = a_i$, we get 
\begin{align*}
\{{{v_i}\over{s-a_i}} {{P'_u}\over{P_u}}(a_i), v(t) {{P'_u}\over{P_u}}(t)\}
& = {1\over{s-a_i}} {{P'_u}\over{\Pi'}}(a_i) v(t) \pa_t \big( 
{2\over{a_i - t}} {{\Pi}\over{P_u}}(t)\big) 
\\& 
+ {{v_i}\over{s-a_i}} {{P'_u}\over{P_u}}(t)
\big( {2\over{a_i - t}} {{\Pi'}\over{P_u}}(a_i) {{P_u}\over{\Pi}}(t)
+ {2\over{(a_i - t)^2}}\big) 
\end{align*}
and 
$$
\{{{v_i}\over{s-a_i}} {{P'_u}\over{P_u}}(a_i), {{v_j}\over{t-a_j}}
{{P'_u}\over{P_u}}(a_j)\}
= {2\over{a_i - a_j}} {1\over{s-a_i}} {1\over{t-a_j}}
\big( v_j {{P'_u}\over{\Pi'}}(a_i) {{\Pi'}\over{P_u}}(a_j)
+ v_i {{P'_u}\over{\Pi'}}(a_j) {{\Pi'}\over{P_u}}(a_i)
\big) 
$$ 
for $i\neq j$ and $= 0$ for $i = j$. 

Since $\Pi'(a_i) = P_u(a_i) / u_i$, we get 
\begin{align*}
\{W(s),W(t)\} & = 
\Big( v(t){{P'_u}\over{P_u}}(s) \pa_t \Big( {2\over{s-t}} \big( 
{{\Pi(t)}\over{\Pi(s)}} {{P_u(s)}\over{P_u(t)}} - 1 \big) \Big)  
\\ & 
- v(t) \sum_i {1\over{s-a_i}} u_i {{P'_u}\over{P_u}}(a_i)
\pa_t \big( {2\over{a_i - t}} {\Pi\over{P_u}}(t) \big)
\\ & 
- {{P'_u}\over{P_u}}(t) \sum_i {{v_i}\over{s-a_i}}
\big( {2\over{a_i - t}} {1\over{u_i}} {{P_u}\over{\Pi}}(t)
+ {2\over{(a_i-t)^2}}\big) \Big) - (s\leftrightarrow t)
\\ & 
- \sum_{i,j\¦ i\neq j} {2\over{a_i - a_j}} {1\over{s-a_i}} {1\over{t-a_j}}
\Big( v_j {{u_i}\over{u_j}} {{P'_u}\over{P_u}}(a_i)
+ v_i {{u_j}\over{u_i}} {{P'_u}\over{P_u}}(a_j)\Big) .
\end{align*}
The second term is rewritten as 
$$
v(t) {{P'_u}\over{\Pi}}(s) \pa_s \big( {{-2}\over{s-t}} {\Pi\over{P_u}}(t)\big) 
+ v(t) \pa_t \big( {2\over{s-t}} {{P'_u}\over{P_u}}(t)\big) . 
$$
After simplifications, we get 
\begin{align*}
& \{W(s),W(t)\} =
{2\over{s-t}} v(t) \big({{P'_u}\over{P_u}} \big)'(t)
\\& 
+ {{P'_u}\over{P_u}}(t) \big( {2\over{(s-t)^2}}(v(s)+v(t))
- \sum_i {{2v_i}\over{(s-a_i)(t-a_i)^2}} \big)  
\\ & 
+ {2\over{s-t}} v(s) \big({{P'_u}\over{P_u}}\big)'(s)
+ {{P'_u}\over{P_u}}(s) \big( - {2\over{(s-t)^2}}(v(s)+v(t))
+ \sum_i {{2v_i}\over{(s-a_i)^2(t-a_i)}} \big)  
\\ & 
+ \big( {{P'_u}\over{\Pi}}(t) - {{P'_u}\over{\Pi}}(s) \big) 
\sum_i {{2v_i / u_i}\over{(s-a_i)(t-a_i)}} 
\\ & 
- \sum_{i,j\¦ i\neq j} {2\over{a_i - a_j}} {1\over{s-a_i}} {1\over{t-a_j}}
\Big( v_j {{u_i}\over{u_j}} {{P'_u}\over{P_u}}(a_i)
+ v_i {{u_j}\over{u_i}} {{P'_u}\over{P_u}}(a_j)\Big) .
\end{align*}

On the other hand, the right hand side of (\ref{PB:2bis}) is equal to 
\begin{align*}
& {2\over{s-t}} v(t) \big({{P'_u}\over{P_u}}\big)'(t)
+ {{P'_u}\over{P_u}}(t) \big( {2\over{s-t}} v'(t) + {4\over{(s-t)^2}} v(t)\big) 
\\ & 
+ {2\over{s-t}} v(s) \big({{P'_u}\over{P_u}}\big)'(s)
+ {{P'_u}\over{P_u}}(s) \big( {2\over{s-t}} v'(s) - {4\over{(s-t)^2}} v(s)\big) 
\\ & 
+ {2\over{s-t}} \sum_i v_i {{P'_u}\over{P_u}}(a_i) \big( {1\over{(s-a_i)^2}}
+ {1\over{(t-a_i)^2}}\big) 
+ {4\over{(s-t)^2}} \sum_i v_i {{P'_u}\over{P_u}}(a_i) 
\big( {1\over{s-a_i}} - {1\over{t-a_i}}\big) . 
\end{align*}

Identity (\ref{PB:2bis}) then follows from the identities 
$$
{2\over{s-t}} v'(t) + {4\over{(s-t)^2}} v(t) = 
{2\over{(s-t)^2}} (v(t)+ v(s)) - \sum_i {{2v_i}\over{(s-a_i)(t-a_i)^2}}
$$
and 
\begin{align*}
& \big( {{P'_u}\over{\Pi}}(t) - {{P'_u}\over{\Pi}}(s) \big) 
\sum_i {{2v_i / u_i}\over{(s-a_i)(t-a_i)}} 
\\& - \sum_{i,j\¦i\neq j}
{2\over{a_i - a_j}} {1\over{s-a_i}} {1\over{t-a_j}}
\Big( v_j {{u_i}\over{u_j}} {{P'_u}\over{P_u}}(a_i)
+ v_i {{u_j}\over{u_i}} {{P'_u}\over{P_u}}(a_j)\Big) 
\\& 
= {2\over{s-t}} \sum_i v_i {{P'_u}\over{P_u}}(a_i)
\big( {1\over{s-a_i}} - {1\over{t-a_i}}\big)^2 .
\end{align*}
The first identity is immmediate, and the second identity follows from 
the identity 
$$
{{P'_u}\over{\Pi}}(t) = \sum_i {1\over{t-a_i}} u_i {{P'_u}\over{P_u}}(a_i). 
$$
\hfill \qed \medskip 

This ends the proof of Theorem \ref{thm:iso:class}. 
\hfill \qed \medskip 

\subsection{}

The isomorphic ringed Lie algebras $\cX \simeq \cU$ may be described
as follows. 
Let us define $P_{\ol u} \in K[X]$ as 
$$
P_{\ol u}(X) = \prod_{i = 1}^{g+3} (X - a_i) \cdot 
\sum_{i = 1}^{g+3} {{\ol u_i}\over{X - a_i}} .  
$$
Set $\cA^0 = K$, $\cA^1 = K[X] / (P_{\ol u})$. 
Since $P_{\ol u}$ is a monic polynomial of degree $g$, 
$\cA^1$ is a $g$-dimensional $K$-vector space with basis 
$(1,X,\ldots,X^{g-1})$. Let us set 
$$
\cA = \cA^0 \oplus \cA^1; 
$$
for $\phi \in K[X]$, we denote by $V[\phi]$ the element of 
$(0,\on{class\ of\ }\phi)$ of $\cA$. 
We define a bracket 
$$
\{ , \} : \wedge^2(\cA) \to \cA
$$
of degree $-1$ as follows: 
\begin{prop} \label{prop:A}
1) For each $\phi \in K[X]$, there is a unique derivation 
$\delta(\phi)$ of $K$, such that for any $i = 1,\ldots,g+3$, 
$$
\delta(\phi)(\ol u_i) =  \Big( 
\phi(a_i) {{P'_{\ol u}}\over{P_{\ol u}}}(a_i)
- \big( \phi(X) {{P'_{\ol u}}\over{P_{\ol u}}}(X) \big)_{\geq 0 \¦ X = a_i} 
 \Big) \ol u_i ; 
$$
here the argument of $(\ )_{\geq 0}$ is a Laurent formal series in 
$K((X^{-1}))$, and the index $(\ )_{\geq 0\¦ X = a_i}$ means projection on 
$K[X]$ parallel to $X^{-1} K[[X^{-1}]]$, followed by evaluation of a polynomial 
at $a_i$. Then $\delta$ factors through 
a linear map $\delta : K[X] / (P_{\ol u}) \to \on{Der}(K)$. For 
$\phi\in K[X]$, $u\in K$, we set 
\begin{equation} \label{part:A}
\{V[\phi],u\} = \delta(\phi)(u). 
\end{equation}

2) There is a unique ringed
Lie algebra map $\wedge^2(\cA) \to \cA$, extending (\ref{part:A}) and such 
that for any $\phi,\psi\in \CC[X]$, we have  
$$
\{V[\phi],V[\psi]\} = V[2(\phi' \psi - \phi \psi')(X)]. 
$$
\end{prop}

The isomorphism $\cA \to \cX$ is given by 
$$
\ol u_i \mapsto {{\prod_\al (a_i - x_\al)}\over{\prod_{j\¦\neq i} (a_i - a_j)}} , 
$$
$$
V[\phi] \mapsto Y[\phi], 
$$
and the isomorphism $\cA \to \cU$ is given by 
$$
\ol u_i \mapsto \on{class\ of\ }u_i / \sum_i a_i^2 u_i, 
$$
$$
V[\phi] \mapsto W[\phi]. 
$$

\section{An isomorphism of Poisson fields}
\label{sect:5}

\begin{thm} \label{thm:iso}
The field of rational functions on $\cN^{g+3}/\!/\Gamma$ is 
$\Frac(S^\bullet_{\cU^0}(\cU^1))$; the field of rational functions
on $(T^*(\PP^1))^{(g)}$ is $\Frac(S^\bullet_{\cX^0}(\cX^1))$. 
The isomorphism of ringed Lie algebras $\cX \to \cU$
(see Theorem \ref{thm:iso:class}) induces an isomorphism of gra\-ded Poisson 
algebras $S^\bullet_{\cX^0}(\cX^1) \to S^\bullet_{\cU^0}(\cU^1)$, 
and therefore an isomor\-phism of Pois\-son fields
$$
\Frac(S^\bullet_{\cX^0}(\cX^1)) \to \Frac(S^\bullet_{\cU^0}(\cU^1)), 
$$
i.e., a birational Poisson morphism 
$\cN^{g+3}/\!/\Gamma \to (T^*(\PP^1))^{(g)}$. 
\end{thm}

\begin{prop} \label{prop:iso}
The Poisson-commuting families $(H_1^\BM,\ldots, H_g^\BM)$
and $(\ul H_1^\Hitch,$ $\ldots, \ul H_{g+3}^\Hitch)$ (see Section \ref{sect:1})
actually belong to $S^2_{\cX^0}(\cX^1) \subset S^\bullet_{\cX^0}(\cX^1)$
and to $S^2_{\cU^0}(\cU^1)$ $\subset S^\bullet_{\cU^0}(\cU^1)$. 

The isomorphism of ringed Lie algebras $\cX \to \cU$ induces an 
isomorphism of graded Poisson algebras $S^\bullet_{\cX^0}(\cX^1)
\to S^\bullet_{\cU^0}(\cU^1)$, which takes $H_i^\BM$ to 
$$
(-1)^{g-i} \sum_{j=1}^{g+3} \ul H_j^\Hitch
\sum_{\stackrel{
\stackrel{(j_1,\ldots,j_{g+2-i}) \¦}{ j_1,\ldots,j_{g+2-i} \in 
\{1,\ldots,\check j,\ldots,g+3\}}} 
{j_1 < \ldots < j_{g+2-i}} }
a_{j_1}\ldots a_{j_{g+2-i}} . 
$$
\end{prop}

{\em Proof.} One can show that 
$S^\bullet_{\cX^0}(\cX^1)$ is the subalgebra 
$$
\CC(x_1,\ldots,x_g)[y_1,\ldots,y_g]^{\SG_g}
\subset \CC(x_1,\ldots,y_g)^{\SG_g}. 
$$
On the other hand, $S^\bullet_{\cU^0}(\cU^1)$ is the reduction of 
$$
\cO = \CC(u_1,\ldots,u_{g+3})_{\on{Span}(\ell_0,\ell_1)}[v_1,\ldots,v_{g+3}]
$$
with respect to the ideal $I = (\sum_i v_i, \sum_i u_i, \sum_i a_i u_i)$;
here the reduction is defined as 
$$
\{f\in \cO \¦ \{f,I\} \subset I\} / I. 
$$ 
Set 
$$
\wt H_i^\Hitch = \sum_{j\¦ j\neq i} {1\over{a_j - a_i}}
{1\over {u_i u_j}}(u_i v_j - u_j v_i)^2 ; 
$$
then $\wt H^\Hitch_i\in \cO$, and $\{\wt H^\Hitch_i,I\} \subset I$, 
therefore $\wt H^\Hitch_i$ defines elements in the reduced algebra; 
these elements are the $\ul H_i^\Hitch$, $i = 1,\ldots,g+3$. 
\hfill \qed \medskip 


\section{Quantization of the isomorphism of Poisson algebras}
\label{sect:quant}

We summarize Theorem \ref{thm:iso} and Proposition \ref{prop:iso} 
by saying that there is an isomorphism of integrable systems between 
$$
\big(S^\bullet_{\cX^0}(\cX^1),(H_1^\cX,\ldots,H_g^\cX)\big) \; \on{and} \; 
\big(S^\bullet_{\cU^0}(\cU^1),(H_1^\cU,\ldots,H_{g+3}^\cU)\big). 
$$
By a quantization of this isomorphism, we understand the following 
data: 

(1) explicit quantizations $A_\cX$ and $A_\cU$ of the Poisson algebras
$S^\bullet_{\cX^0}(\cX^1)$ and $S^\bullet_{\cU^0}(\cU^1)$, and an 
isomorphism $i : A_\cX \to A_\cU$ quantizing the isomorphism 
$S^\bullet_{\cX^0}(\cX^1) \to S^\bullet_{\cU^0}(\cU^1)$; 

(2) commuting families $(\wh H_1^\cX,\ldots,\wh H_g^\cX)$ of $A_\cX$
and $(\wh H_1^\cU,\ldots,\wh H_{g+3}^\cU)$ of $A_\cU$, quantizing 
$(H_1^\cX,\ldots,H_g^\cX)$ and $(H_1^\cU,\ldots,H_{g+3}^\cU)$, which are
mapped to each other by $i$. 

We will construct a family of quantizations of this isomorphism, associated with 
a collection $(C_1,\ldots,C_{g+3})$ of scalars (values of the Casimir elements 
at the marked points). 

\subsection{Quantization of $S^\bullet_{\cX^0}(\cX^1)$}

We have 
$$
S^\bullet_{\cX^0}(\cX^1) = \big( 
\CC(x_1,\ldots,x_g)[y_1,\ldots,y_g]\big)^{\SG_g}. 
$$

Let $\CC(\wh x_1,\ldots,\wh x_g)[\wh y_1,\ldots,\wh y_g]$ be the semidirect 
product of the commutative field $\CC(\wh x_1,\ldots,\wh x_g)$ by the abelian 
Lie algebra $\oplus_{i=1}^g \CC \wh y_i$, acting by the derivations 
$\wh y_i \mapsto (f \mapsto 2\Pi(\wh x_i) \pa_{\wh x_i} f)$. 
So we have the relations 
$$
[\wh y_i,\wh y_j] = 0, \; 
[\wh y_i, f(\wh x_1,\ldots,\wh x_g)] = 
2\Pi(\wh x_i) {{\pa f}\over{\pa \wh x_i}}
(\wh x_1,\ldots,\wh x_g). 
$$

Then $\SG_g$ acts on $\CC(\wh x_1,\ldots,\wh x_g)[\wh y_1,\ldots,\wh y_g]$ 
by permuting the pairs $(\wh x_i,\wh y_i)$. 

\begin{prop}
Set $A_\cX = \big( \CC(\wh x_1,\ldots,\wh x_g)[\wh y_1,\ldots,\wh y_g] 
\big)^{\SG_g}$. Then $A_\cX$ is a quantization of $S^\bullet_{\cX^0}(\cX^1)$. 

For $\phi \in \CC(x_1,\ldots,x_g)^{\SG_g}$ and $\psi\in 
\CC(x_1,\ldots,x_g)^{\SG_{g-1}}$ (where $\SG_{g-1}$ permutes the $g-1$ last
variables), set 
$$
\wh X[\phi] = \phi(\wh x_1,\ldots,\wh x_g), \; 
\wh Y[\psi] = \sum_\al \psi(\wh x_\al,\wh x_1,\ldots,\wh x_{\al-1},
\wh x_{\al +1},\ldots,\wh x_g) \wh y_\al. 
$$ 

Then the $\wh X[\phi]$ and $\wh Y[\psi]$ generate $A_\cX$. 
Relations between these generators are: 
\begin{equation} \label{rel:linear}
\phi \mapsto \wh X[\phi] \; \on{and}\; 
\psi \mapsto \wh Y[\psi]\, \on{are\ linear,}
\end{equation}
\begin{equation} \label{rel:X:Y:1}
\wh X[\phi \phi'] = \wh X[\phi] \wh X[\phi], \; 
\wh Y[\phi \psi] = \wh X[\phi] \wh Y[\psi], 
\end{equation}
\begin{equation} \label{rel:X:Y:2}
[\wh X[\phi],\wh Y[\psi]] = 
-2 \wh X[\sum_\al \Pi(x_\al) \psi(x_\al,x_1,\ldots,\check{x_\al},\ldots,x_g)
(\pa_\al \phi)(x_1,\ldots,x_g)] , 
\end{equation}
\begin{equation} \label{rel:X:Y:3}
[\wh Y[\psi],\wh Y[\psi']] = 
2 \wh Y[\sum_\al \Pi(x_\al) \psi(x_\al,x_1,\ldots,\check{x_\al},\ldots,x_g)
(\pa_\al \psi')(x_1,\ldots,x_g) - (\psi \leftrightarrow \psi') ] , 
\end{equation}
for any $\phi,\phi'\in \CC(x_1,\ldots,x_g)^{\SG_g}$
and $\psi,\psi'\in \CC(x_1,\ldots,x_g)^{\SG_{g-1}}$. 

$A_{\cX}$ is isomorphic to the algebra generated by the $\wh X[\phi]$
and $\wh Y[\psi]$, subject to relations (\ref{rel:linear}), (\ref{rel:X:Y:1}),   
(\ref{rel:X:Y:2}) and (\ref{rel:X:Y:3}). 
\end{prop}

\subsection{Quantization of Poisson algebras associated with vector spaces}

Let $W\subset V$ be an inclusion of finite-dimensional 
vector spaces. Let us set 
$\cO = \CC(V)_W$ and $I = (W) \subset \cO$. We have a Lie algebra 
morphism $V^* \to \Der(\cO)$, $\xi \mapsto \pa_\xi$, where $V^*$
is equipped with the zero Lie bracket. Let $A = \cO[V^*]$ be the 
corresponding crossed product algebra. Then $A$ is a localization of 
the Weyl algebra of $V$. The product induces an isomorphism 
$\cO \otimes S^\bullet(V^*) \to A$. 

Let us set 
$$
N = \{f\in A \¦ f W \subset WA\}. 
$$
Then $N$ is a subalgebra of $A$; the product of $A$ induces an 
isomorphism $(W) \otimes S^\bullet(V^*) + \cO \otimes S^\bullet(W^\perp)
\to N$. 
Moreover, $WA \subset N$ is a two-sided ideal of $N$, and the quotient 
$A_{\on{red}} = N / WA$ is isomorphic to 
$$
A_{\on{red}} = \CC(V/W)[(V/W)^*],  
$$
a localization of the Weyl algebra of $V/W$.  

Let now $\Phi = \sum_i x_i \otimes \xi^i \in V \otimes V^*$ be the canonical 
element. Define a grading on a part of $\CC(V)[V^*]$ by giving degree $1$ 
to the elements of $V\subset \CC(V)$ and degree $-1$ to the elements of 
$V^* \subset \CC(V)[V^*]$. We denote by $\CC(V)[V^*]^0$ the degree zreo 
part of $\CC(V)[V^*]$. Then if $\CC(V)^\al$ is the degree $\al$ part of 
$\CC(V)$, then the product induces a linear isomorphism 
$$
\bigoplus_{\al\geq 0} \CC(V)^{\al} \otimes S^\al(V^*) \to \CC(V)[V^*]^0. 
$$
Then $\Phi$ is a central element in $\CC(V)[V^*]^0$. 
The quantization of the Poisson algebra associated to 
(\ref{G:V:0}) is 
$$
\CC(V)[V^*]^0  / (\Phi). 
$$

Applying this reduction to $A_{\on{red}}$, we obtain the following quantization 
$A'_{\on{red}}$ of the Poisson algebra associated with (\ref{reduced:GLA}): set 
$N^0 = \{f\in A^0 \¦ W f \subset A^0 W\}$, then $N^0$ contains 
$A^{-1} W$ as a 2-sided ideal, and $\Phi$ as a central element. Then  
$$
A'_{\on{red}} = N^0 / \big( A^{-1} W + N^0 \Phi\big).  
$$

\subsection{Quantization of $S^\bullet_{\cU^0}(\cU^1)$}

A quantization of $S^\bullet_{\cU^0}(\cU^1)$ is then the algebra 
$A'_{\on{red}}$ associated to $V = \oplus_{i=1}^{g+3} \CC f_i$
and $W = \on{Span}(\wh \ell_0,\wh \ell_1)$, where $\wh \ell_0 
= \sum_i f_i$ and 
$\wh \ell_1 = \sum_i a_i f_i$. 

We denote by $A,A_{\on{red}}$ and $A_\cU = A'_{\on{red}}$ the resulting algebras. 
The relations 
$$
[f_i^*,f_j] = \delta_{i,j}, \; [f_i,f_j] = 
[f_i^*,f_j^*] = 0 
$$
hold in $A$. We set $h_i = -2 f_i f_i^* \in A$. 

We also set $\ol f_i = f_i / (\sum_i a_i^2 f_i)$. Then 
we have $\ol f_i \in N^0$. We also denote by $f_i,\ol f_i$ the images of 
$f_i,\ol f_i$ in $A'_{\on{red}}$. Then we have the relations 
$$
\sum_i \ol f_i = \sum_i a_i \ol f_i = 0, \; 
\sum_i a_i^2 \ol f_i = 1 
$$ 
in $A'_{\on{red}}$. We therefore have an injection $K \hookrightarrow 
A'_{\on{red}}$. 

Elements of degree $\leq 1$ of $A'_{\on{red}}$ are 
the classes of the $\sum_i \la_i f_i^*$, where 
$(\la_1,\ldots,\la_{g+3}) \in K^{g+3}$ is such that 
$\sum_i \la_i = \sum_i a_i \la_i = 0$. 

\medskip 

$A_\cU$ may also be considered as the reduction of 
$$
A = \CC(f_1,\ldots,f_{g+3})_{\on{Span}(\wh \ell_0,\wh \ell_1)}
[h_1,\ldots,h_{g+3}]
$$ 
with respect to the left ideal 
$$
\wh I = A (\sum_i h_i - 2) + A \wh \ell_0 + A \wh \ell_1. 
$$
Then 
\begin{equation} \label{hecke}
A_\cU = \wh N / \wh I,  
\end{equation}
where $\wh N = \{f\in A \¦ I f \subset I\}$.

\subsection{Quantization of the isomorphism $S^\bullet_{\cX^0}(\cX^1)
\to S^\bullet_{\cU^0}(\cU^1)$}

We set 
$$
P_f(X) = \prod_{i} (X - a_i) \cdot \sum_{i} {{f_i}\over{X-a_i}}
\in A_{\on{red}}[X]. 
$$

\begin{thm} \label{thm:quant:GNR}
1) For any $\phi \in \CC[X]$, the element 
$$
\wh W[\phi] = \on{Res}_{t = \infty} \Big( \sum_i 
\big( {{P'_{f}}\over{P_{f}}}(t) - {{P'_{f}}\over{P_{f}}}(a_i)\big) 
{{h_i}\over{t-a_i}} \phi(t)dt\Big) 
$$
belongs to $N^0$; we also denote by $\wh W[\phi]$ its image in $A'_{\on{red}} = 
A_{\cU}$. 

2) There is a unique isomorphism $\wh \La : A_\cX \to A_\cU$, taking 
$$
\sum_{1\leq \al_1 < \ldots < \al_k \leq g} \wh x_{\al_1} \cdots \wh x_{\al_k}
\mapsto  \sum_i \ol f_i 
\sum_{ \stackrel{j_1 < \ldots < j_{k+2}}{j_1,\ldots,j_{k+2} \neq i} } 
a_{j_1} \cdots a_{j_{k+2}}
$$
and 
$$
\wh Y[\phi] = \sum_\al (\phi/\Pi)(\wh x_\al) \wh y_\al
\mapsto 
\wh W[\phi] = \on{Res}_{t = \infty} \Big( \sum_i 
\big( {{P'_{f}}\over{P_{f}}}(t) - {{P'_{f}}\over{P_{f}}}(a_i)\big) 
{{h_i}\over{t-a_i}} \phi(t)dt\Big).  
$$
\end{thm}

{\em Proof.} Set 
$$
\wh W(X) = \sum_i \big( {{P'_{f}}\over{P_{f}}}(X) - 
{{P'_{f}}\over{P_{f}}}(a_i)\big) {{h_i}\over{X-a_i}}
\in X^{-1} A_{\on{red}}[[X^{-1}]].   
$$
Then
$$
\wh W[\phi] = \on{Res}_{t=\infty}(\wh W(t)\phi(t)dt). 
$$
We have 
$$
[\wh \ell_0,\wh W(X)] = 2 {{P'_{f}(X)}\over{\Pi(X)}}
- 2 \sum_i {{f_i}\over{X-a_i}} {{P'_{f}}\over{P_{f}}}(a_i) = 0
$$
by the same argument as above. 
We have also 
$$
[\wh \ell_1,\wh W(X)] = - 2 \wh \ell_0 {{P'_{f}}\over{P_f}}(X)
+2 \sum_i f_i {{P'_{f}}\over{P_{f}}}(a_i). 
$$
As before, we derive from this $\sum_i f_i (P'_{f}/P_{f})(a_i) = 0$, 
therefore $\wh W(X)\in X^{-1} N[[X^{-1}]]$. 

On the other hand, $\wh W(X)$ commutes with $\sum_{i=1}^{g+3} h_i$, 
so it belongs to $X^{-1} N^0[[X^{-1}]]$. This proves 1). 

Let us prove 2). Let us set 
$$
\wh X(t) = \sum_\al {1\over{t-\wh x_\al}} \in t^{-1} A_\cX[[t^{-1}]], \; 
\wh Y(t) = \sum_\al {{\Pi(\wh x_\al)^{-1}}\over{t - \wh x_\al}} \wh y_\al 
\in t^{-1} A_\cX[[t^{-1}]], 
$$
and 
$$
\wh U(t) = {{P'_{f}(t)}\over{P_{f}(t)}} \in t^{-1} A_\cU[[t^{-1}]], \; 
\wh W(t) = \sum_i \big( {{P'_{f}}\over{P_{f}}}(t)
- {{P'_{f}}\over{P_{f}}}(a_i)\big) {{h_i}\over{t-a_i}}
\in t^{-1} A_\cU[[t^{-1}]].  
$$

Then $\wh X(t), \wh Y(t)$ satisfy the relations 
$$
[\wh X(t),\wh X(s)] = 0, 
$$
$$
[\wh Y(t),\wh X(s)] = 2 \pa_s \big( {{\wh X(s)-\wh X(t)}\over{s-t}}\big) , 
$$
$$
[\wh Y(t),\wh Y(s)] = {2\over{s-t}} \big( \pa_s\wh Y(s) + \pa_t\wh Y(t)
- {2\over{s-t}} (\wh Y(s) - \wh Y(t))  \big).  
$$
A family of relations for $A_\cX$ can be extracted from these relations. 
We have $\wh \La(\wh X(t)) = \wh U(t)$, $\wh \La(\wh Y(t)) = \wh W(t)$. 

In the same way as the proof of Theorem \ref{thm:iso:class}, one shows 
$$
[\wh U(t),\wh U(s)] = 0, 
$$
$$
[\wh W(t),\wh U(s)] = 2 \pa_s \big( {{\wh U(s)-\wh U(t)}\over{s-t}}\big) , 
$$
$$
[\wh W(t),\wh W(s)] = {2\over{s-t}} \big( \pa_s\wh W(s) + \pa_t\wh W(t)
- {2\over{s-t}} (\wh W(s) - \wh W(t))  \big).  
$$
This implies that $\wh\La : A_\cX \to A_\cU$ is an algebra morphism. Since 
its associated graded is the isomorphism $\La$, $\wh\La$ is also an 
isomorphism. 
\hfill \qed \medskip 

\begin{remark} \label{rem:hecke}
$A_\cU$ contains the Hecke algebra $B/J$, where 
$B \subset \otimes_{i=1}^N \big(U(\SL_2) / (C - C_i)\big)$
is $\{f\¦ Jf \subset J\}$, and $J$ is the left ideal generated by 
$\sum_i h_i-2$, $\sum_i f_i$ and $\sum_i a_i f_i$. $B/J$ may be viewed
as the quantization of the algebra of regular functions on 
$\cN^{g+3}//\Gamma$. So the isomorphism $A_\cU \to A_\cX$ restricts
to an isomorphism of quantized regular function algebras. 
\end{remark}

\subsection{Quantization of the Hitchin hamiltonians}

Let $C_1,\ldots,C_{g+3}$ be scalars. Let us set 
$$
e_i = - {1\over {2 f_i}} C_i - {1\over{4f_i}} h_i^2 - {1\over{2 f_i}} h_i. 
$$
Then there is a unique algebra morphism 
$$
\bigotimes_{i=1}^{g+3} \big( U(\SL_2) / 
(ef + fe+ {1\over 2} h^2 - C_i) \big) \to A_\cU, 
$$
taking each $x^{(i)}$ to $x_i$, $x\in \{e,f,h\}$, $i\in \{1,\ldots,g+3\}$. 

Let us set 
$$
\wh H_i^{\on{Hitch}} = \sum_{j\¦j\neq i}
{{e_i f_j + f_i e_j + {1\over 2} h_i h_j}\over{a_i - a_j}}. 
$$
Then if we set 
$$
x(X) = \sum_i {{x_i}\over{X-a_i}}, 
$$
for $x\in\{e,f,h\}$, we have the identity 
$$
\sum_{i=1}^{g+3} {{\wh H_i^{\on{Hitch}}}\over{X-a_i}} + \sum_{i=1}^{g+3} 
{{C_i/2}\over{(X-a_i)^2}} = {1\over 2} \big( e(X)f(X) + f(X)e(X) 
+ {1\over 2} h(X)^2 \big) . 
$$

\begin{prop} \label{prop:hecke} \label{prop:QIS}
The $\wh H_i^{\on{Hitch}}$ belong to $\wh N = \{ f\in A \¦ I f \subset I\}$. 
Their classes in $A_\cU$, also denoted $\wh H_i^{\on{Hitch}}$, form a 
commutative family, quantizing the $H_i^{\on{Hitch}}$. Moreover, we have 
\begin{equation} \label{ids:H:Hitch}
\sum_i \wh H_i^{\on{Hitch}} = 0, \; 
\sum_i a_i \wh H_i^{\on{Hitch}} = - {1\over 2} \sum_i C_i, \; 
\sum_i a_i^2 \wh H_i^{\on{Hitch}} = - \sum_i a_i C_i \; 
\end{equation}
(equalities in $A_\cU$).  
\end{prop}

{\em Proof.}
Let us set 
$$
\wh H(X) = e(X)f(X) + f(X)e(X) + {1\over 2} h(X)^2.  
$$
We have 
$$
(\sum_i h_i - 2)\wh H(X) = 
\wh H(X) (\sum_i h_i - 2), 
$$
$$
(\sum_i f_i)\wh H(X) = 
\wh H(X) (\sum_i f_i), 
$$
and 
$$
(\sum_i a_i f_i) \wh H(X) 
= \wh H(X) (\sum_i a_i f_i) - 2 f(X) (\sum_i h_i - 2) 
+ 2 h(X) (\sum_i f_i). 
$$
This proves that the $\wh H_i^{\on{Hitch}}$ belong to $\wh N$. 

On the other hand, we have 
$$
\wh H(X) = e(X)f(X) + {1\over 2} h(X)^2 - h'(X), 
$$
whose expansion in $X^{-1}$ gives 
$$
\sum_i \wh H_i^{\on{Hitch}} = 0, \;  
\sum_i a_i \wh H_i^{\on{Hitch}} + \sum_i C_i / 2
= (\sum_i e_i)(\sum_i f_i) + {1\over 4}(\sum_i h_i)(\sum_i h_i - 2), \;  
$$
and 
$$
\sum_i a_i^2 \wh H_i^{\on{Hitch}} + \sum_i a_i C_i 
= (\sum_i a_i e_i)(\sum_i f_i) + 
(\sum_i e_i)(\sum_i a_i f_i) +{1\over 2}(\sum_i a_i h_i)(\sum_i h_i - 2). 
$$
The right sides of these identities all belong to $\wh I$. This 
implies the identities (\ref{ids:H:Hitch}). 
\hfill \qed \medskip 

\subsection{Quantization of the Beauville-Mukai hamiltonians}

Let $\wh T$ be an element of the algebra $A = \CC(\wh x)[\wh y]$, 
where the relations are 
$$
[\wh y,f(\wh x)] = 2\Pi(\wh x) (\pa_{\wh x}f)(\wh x).
$$ 
The Beauville-Mukai hamiltonians corresponding to the 
family $(1,\wh x,\ldots, \wh x^{g-1},\wh T)$ are the 
$$
\wh H_i^{\on{BM}} = \big( \prod_{j>1} (\wh x_j - \wh x_1)\big)^{-1} 
\sum_{1<j_1<\ldots<j_{g-i}} \wh x_{j_1} \cdots \wh x_{j_{g-i}} \wh T_1
+ \on{\ cyclic\ permutations}.
$$
We set $f_i = f^{(i)} = 1^{\otimes i-1} \otimes f \otimes 1^{\otimes g-i}\in 
A^{\otimes g}$. 

They satisfy the identity 
$$
\sum_{i=1}^g {{\prod_{j\¦j\neq i} (X-\wh x_j)}
\over {\prod_{j\¦j\neq i} (\wh x_i-\wh x_j)}} \wh T_i = 
\sum_{i=0}^{g-1} (-1)^i X^i \wh H_{i+1}^{\on{BM}}. 
$$ 
in $A^{\otimes g}[X]$. 

We will set 
$$
\wh T = (\Pi(\wh x_\al)\wh y_\al)^2 + \sum_{i=1}^{g+3} {{-C_i/2}\over{(\wh x - a_i)^2}}. 
$$
Then $(\wh H_i^{\on{BM}})_{i=1,\ldots,g}$ is a commuting family of 
quantum hamiltonians, quantizing $(H_i^{\on{BM}})_{i=1,\ldots,g}$.  

\subsection{Isomorphism of quantum integrable systems}

\begin{thm}
Assume that $(C_i)_{i=1,\ldots,g+3}$ satisfies $\sum_i C_i = \sum_i a_i C_i$, 
then the families $(\wh H_i^{\on{Hitch}})_{i=1,\ldots,g+3}$ 
and $(\wh H_i^{\on{BM}})_{i=1,\ldots,g}$ are mapped to each other by the isomorphism 
$\wh\La$.  More precisely, we have the identities 
$$
\sum_{i=0}^{g-1} (-1)^i X^i \wh\La(\wh H_{i+1}^{\on{BM}}) = 
\big( \sum_i {{\wh H_i^{\on{Hitch}}}\over{X-a_i}}\big) \prod_{i=1}^{g+3} (X-a_i),  
$$
that is 
$$
\wh H_i^{\on{Hitch}} = \big( \prod_{j\¦j\neq i} (a_i - a_j)\big)^{-1}
\sum_{k=0}^{g-1} (-1)^k (a_i)^k \wh\La(\wh H_{k+1}^{\on{BM}})
$$
and
$$
\wh H_k^{\on{BM}} = (-1)^g \sum_i \wh\La^{-1}(\wh H_i^{\on{Hitch}})
\sum_{\stackrel{j_1< \ldots < j_{g+3-k}}{j_1,\ldots,j_{g+3-k}\neq i}}
a_{j_1} \cdots a_{j_{g+3-i}}. 
$$
\end{thm}

{\em Proof.} Let $A^{\on{big}}$ be the algebra with generators 
$$
\wh x_\al,(\wh x_\al - a_i)^{-1}, \al = 1,\ldots,g, \: 
e_i,h_i,f_i, i = 1,\ldots,g+3,  
$$
and the following relations: 

(1) there is a morphism $\bigotimes_{i=1}^{g+3} (U(\SL_2) / (C-C_i)) \to 
A^{\on{big}}$, taking $x^{(i)}$ to $x_i$, $i = 1,\ldots,g+3$, 
$x\in \{e,h,f\}$; 

(2) we have $\sum_i f_i = \sum_i a_i f_i = 0$, and 
$$
(\sum_i a_i^2 f_i) \prod_\al (X - \wh x_\al) = 
\prod_i (X-a_i) \sum_i {{f_i}\over{X-a_i}}.   
$$

Set for any $\al$, 
\begin{equation} \label{Skl:1}
\wh y_\al = - \Pi(\wh x_\al) \sum_i (\wh x_\al - a_i)^{-1} h_i. 
\end{equation}

Let $A_{\ul a}$ be the localization in $\wh x-a_i$ of 
$\CC\langle \wh x,\wh y\rangle / ([\wh y,\wh x] = 2 \Pi(\wh x))$, then 
there is a morphism $(A_{\ul a})^{\otimes g} \to A^{\on{big}}$, such that 
$\wh x^{(\al)} \mapsto \wh x_\al$, $((\wh x-a_i)^{-1})^{(\al)} 
\mapsto (\wh x_\al-a_i)^{-1}$, $\wh y^{(\al)} \mapsto \wh y_\al$.   

The subalgebras $B/J \subset A_\cU$ and $(A_{\ul a})^{\otimes g} 
\subset A_\cX$ are both subalgebras of $A^{\on{big}}$. These inclusions
are compatible with the isomorphism $\wh\La : A_\cX \to A_{\cU}$. 

\begin{lemma} \label{lemma:end}
We have 
$$
\sum_\al \wh\La^{-1}(\sum_\al \phi(\wh x_\al)(\wh x_\al - a_i)^{-1})
= \wh\La^{-1}(- \sum_\al (\phi / \Pi)(\wh x_\al)\wh y_\al) 
$$
\end{lemma}

{\em Proof.} This follows from 
$$
\sum_\al \phi(\wh x_\al)(\wh x_\al - a_i)^{-1} = 
- \on{Res}_{t = \infty} \big( \big({{P'_f}\over{P_f}}(t) 
- {{P'_f}\over{P_f}}(a_i)  \big)  {{\phi(t)}\over{t-a_i}} dt\big) .
$$ 
\hfill \qed \medskip 

We now show: 
\begin{equation} \label{Skl:2}
\sum_i (\wh x_\al - a_i)^{-1} \wh H_i^{\on{Hitch}} = 
(\Pi(\wh x_\al)^{-1}\wh y_\al)^2 
+ \sum_i {{-C_i / 2}\over{(\wh x_\al - a_i)^2}}
\end{equation}
(equality in $A^{\on{big}}$). This will imply the desired identities. 

Set, for $x\in \{e,f,h\}$, $x(\wh x_\al) = \sum_\al (\wh x_\al - a_i)^{-1} x_i$. 
Then $f(\wh x_\al) = 0$, and $h(\wh x_\al) = - \Pi(\wh x_\al)^{-1} \wh y_\al$. 

We follow Sklyanin's computation (see \cite{Skl}): 
\begin{align*}
\sum_i (\wh x_\al - a_i)^{-1} \wh H_i^{\on{Hitch}} & 
= f(\wh x_\al) e(\wh x_\al) + {1\over 4} h(\wh x_\al)^2
+ \sum_i {{-C_i / 2}\over{(\wh x_\al - a_i)^2}}
\\& + {1\over 2} \sum_i (\wh x_\al-a_i)^{-1} h_i  
+ {1\over 4} \sum_{i,j} (\wh x_\al - a_j)^{-1} 
[(\wh x_\al - a_j)^{-1}, h_j] h_i .  
\end{align*}
Then the first terms vanishes, and the two last terms cancel out due to 
Lemma \ref{lemma:end}. 
\hfill \qed \medskip


\section{Relation with the works of Sklyanin and E. Frenkel} \label{sect:skl}

\subsection{Sklyanin's work}

\subsubsection{Classical case}

In \cite{Skl}, Sklyanin constructs the separation of variables 
for the Gaudin system. In the classical case (a system on $\cN^{g+3}$), 
his map coincides with the map $\GNR$, i.e., it takes the matrix 
$M(X) = \pmatrix A(X) & C(X) \\B(X) & -A(X) \endpmatrix$ to  
the set $\{(x_i,y_i),i = 1,\ldots,g\}$, where 
$\{x_i\}$ is the set of zeroes of $B(X)$ and $y_i = B(x_i)$. 

Let $H^\Hitch(X) = \sum_{i=1}^{g+3} {{H^\Hitch_i}\over{X-a_i}}$
be the generating function of the Hitchin hamiltonians. Set 
$H_\al = H^\Hitch(x_\al)$. 

Sklyanin shows that the $(H_\al)_{\al=1,\ldots,g}$ 
form a "separated system", i.e., there exists a function $H(x,y)$, 
such that $H_\al = H(x_\al,y_\al)$. So the separated system is obtained 
only after mixing the dynamical variables $x_\al$ with the formal 
variable $X$. The Hitchin system is therefore not isomorphic 
with the separated system, but the birational morphism of phase spaces 
$\GNR$ induces a birational morphism between the zero sets of the 
maps $\ul H^\Hitch$ and $(H_1,\ldots,H_g)$. This suggests
that the Beauville-Mukai system should be viewed as an intermediate 
step between the Hitchin system and its separated version. 
This idea was recently exploited in \cite{BT} in the case of other 
integrable systems (the Neumann model and the Toda chain). 

\subsubsection{Quantum case}

The quantum system is formed by the commuting family 
$(\wh H_i^{\Hitch})_{i = 1,\ldots,g+3}$ of $B/J$. He constructs the 
variables $(\wh x_\al)_{\al = 1,\ldots,g}$ and $\wh y_\al$
by (\ref{Skl:1}), and shows the "separation identity" . 
It implies that if $V$ is a module over the quantum algebras, the systems of 
equations $\wh H_i^\Hitch(\psi) = \la_i \psi$ and 
$$
P_{\la_1,\ldots,\la_{g+3}}^{(\al)}(\psi) = 
\big( (\Pi(\wh x_\al)\wh y_\al)^2 
+ \sum_i {{-C_i / 2}\over{(\wh x_\al - a_i)^2}} 
+ {{-\la_i}\over{\wh x_\al - a_i}}\big) (\psi) = 0
$$ 
are equivalent. (\ref{Skl:2}) also implies that the right ideals generated
by $(\wh H_i^{\Hitch}-\la_i)_{i = 1,\ldots,g+3}$ and by 
$(P_{\la_1,\ldots,\la_g}^{(\al)})_{\al = 1,\ldots,g}$ are the same.

\subsection{Frenkel's work}

In \cite{Frenky}, Frenkel shows that the construction of \cite{Skl}
can be translated as the equivalence of the two constructions of a
$\cD$-module over $\Bun_{2,1}(\PP^1,\ul{a})$, related to the geometric Langlands
problem.

One of these $\cD$-modules is based on Radon transformation (see Drinfeld's construction 
in \cite{Dr:AMJ}, and also \cite{Laumon}), while the second is based on the 
Sugawara construction 
(\cite{BD,FFR}). Both $\cD$-modules are quotients of the form  $\cD / 
\sum_i \cD L_i$. The relation with Sklyanin's work is relation (\ref{Skl:2}), 
which shows that two kinds of right ideals are the same. 

\subsection{Relation with the present work}

As in the classical case, the quantum systems defined by 
the operators $(\wh H^\Hitch_i-\la_i)_{i = 1,\ldots,g+3}$
and $(P^{(\la_1,\ldots,\la_{g+3})}_\al)_{\al = 1,\ldots,g}$
are not equivalent, but their left ideals coincide. On the other hand, 
$\wh\GNR$ is an isomorphism of quantum integrable systems, which 
does not involve Radon transformations. Presumably, it can be 
constructed for Gaudin systems in higher rank.

\begin{appendix}
\section{Conformal blocks in degree $\neq 0$}

In this Appendix, we explain the relation of the construction of a quantum 
integrable system in Proposition \ref{prop:QIS} with the usual "conformal
blocks" approach. 

\medskip 

Let $k\geq 0$ be an integer and let $\la_1,\ldots,\la_N$ be a collection of 
integrable weights for $\G = \SL_2$ (here $N = g+3$). To these data  corresponds 
a line bundle $\cL(\la_1,\ldots,\la_N,k)$ over $\Bun_{(\PP^1, \underline{a})}(2,-1)$. 
Its sections may be described as follows. 

Let $V_1,\ldots,V_N$ be the finite-dimensional representations of $\G$ with 
highest weights $\la_1,\ldots,\la_N$. Then 

\begin{lemma}
We have an injection 
\begin{equation} \label{real}
H^0(\Bun_{(\PP^1, \underline{a})}(2,-1), \cL(\la_1,\ldots,\la_N,k))
\hookrightarrow 
\big( (V_1\otimes \cdots \otimes V_N)^* \big)^{\gamma}, 
\end{equation}
where $\gamma \subset U(\G)^{\otimes N}$ is the Lie algebra 
generated by $\sum_{i=1}^N h^{(i)} + k$, $\sum_{i=1}^N f^{(i)}$ 
and $\sum_{i=1}^N a_i f^{(i)}$.  
\end{lemma}

{\em Proof.} The space of sections of $\cL(\la_1,\ldots,\la_N,k)$
may be described as follows. Let $\wh\G = \G((z_\infty)) \oplus \CC K$. 
Let $\G_+ \subset \wh \G$ be the Lie subalgebra 
$$
\G_+  = (e \otimes z_\infty\CC[[z_\infty]]) \oplus 
(h \otimes \CC[[z_\infty]]) \oplus (f \otimes z_\infty^{-1} 
\CC[[z_\infty]]) \oplus  \CC K,   
$$
and let $\chi : \G_+ \to \CC$ be the character such that $\chi(x \otimes 
z_\infty^\al) =  0$ for any $\al$ and $x\in \{e,f\}$, and $\chi(K) = k$. 

Let $(\rho,\VV)$ be the integrable module over $\wh\G$, generated by a
vector $v\in\VV$, such that $\rho(x)(v) = \chi(x)v$ for any $x\in\G_+$. 
We have an inclusion $\G[z] \subset \wh\G$, given by $x\otimes f(z) \mapsto 
x\otimes f(z_\infty^{-1})$. 

Then 
$$
H^0(\Bun_{(\PP^1, \underline{a})}(2,-1), \cL(\la_1,\ldots,\la_N,k))
= \on{Hom}_{\G[z]} \big( V_1(a_1) \otimes \cdots \otimes V_N(a_N)
\otimes \VV,\CC\big),   
$$
where $V(a)$ is the evaluation module over $\G[z]$ corresponding to 
the $\G$-module $V$ and to $a\in \CC$. 

The map (\ref{real}) takes 
$$
\psi \in \on{Hom}_{\G[z]} \big( V_1(a_1) \otimes \cdots \otimes V_N(a_N)
\otimes \VV,\CC\big)   
$$
to the linear form $\wt \psi : v_1\otimes \cdots \otimes v_N \mapsto 
\psi(v_1\otimes \cdots \otimes v_N \otimes v)$. Then 
$$
\wt \psi((\sum_i h^{(i)})(\otimes_{i=1}^N v_i)\otimes v)
= - \psi((\otimes_{i=1}^N v_i)\otimes v) \otimes (h\otimes 1)(v))
= - k \wt\psi((\otimes_{i=1}^N v_i)\otimes v), 
$$
since $(h\otimes 1)v = ([e\otimes z_\infty,f\otimes z_\infty^{-1}]+K)(v) 
=  k v$. In the same way, we have 
$$
\wt \psi((\sum_i f^{(i)})(\otimes_i v_i))
= \wt \psi((\sum_i a_i f^{(i)})(\otimes_i v_i)) = 0. 
$$
\hfill \qed \medskip 

The usual computation shows that the action of the Sugawara tensor 
$T(z_\infty)$ $(dz_\infty)^2$
on $H^0(\Bun_{(\PP^1, \underline{a})}(2,-1), \cL(\la_1,\ldots,\la_N,k))$ is 
the restriction to $\big( (V_1\otimes \cdots \otimes V_N)^* \big)^{\gamma}$ of 
$e(z_\infty)f(z_\infty) + f(z_\infty)e(z_\infty) + {1\over 2} h(z_\infty)^2$, 
where 
$$
x(z_\infty) = \sum_i { {x^{(i)}} \over{ z_\infty^{-1} - a_i} } dz_\infty.
$$ 

When $k = -2$ (critical level), the Sugawara operators are central, so 
this series generates a commutative family of operators on 
$\big( (V_1\otimes \cdots \otimes V_N)^* \big)^{\gamma}$. 

The Hecke algebra of operators acting on 
$\big( (V_1\otimes \cdots \otimes V_N)^* \big)^{\gamma}$ is the quotient
$B/J$ defined in Remark \ref{rem:hecke}. Proposition \ref{prop:hecke} 
shows that the operator $e(z_\infty)f(z_\infty) + f(z_\infty)e(z_\infty) 
+ {1\over 2} h(z_\infty)^2$ belongs to $A_\cU$, which contains the 
Hecke algebra $B/J$. 

\end{appendix}

\subsection*{Acknowledgements}

{\small We would like to thank E. Frenkel, A. Odesskii and M. Olshanetsky 
for discussions on the subject of this work. V.R. would also like to thank 
IRMA (CNRS, Strasbourg) and IHES for hospitality at the time this work was 
done. V.R. was partially supported by grants INTAS 99-1705, RFBR 01-01-00549 
and by the grant for scientific schools RFBR 00-15-96557. }

\end{document}